\documentclass[A4,twoside,11pt,leqno]{article}
\usepackage{amssymb}
\usepackage{amsthm}
\usepackage[tbtags]{amsmath}
\usepackage{doc}
\usepackage{latexsym}
\usepackage{amscd}
\usepackage[frame,cmtip,arrow,matrix,line,graph,curve]{xy}
\usepackage{graphics}
\usepackage{epsfig}
\usepackage{eucal}
\usepackage{mathrsfs}
\usepackage{tikz-cd}
\font\headd=cmr8

\newtheorem{theorem}{Theorem}[section]

\newtheorem{remark}{Remark}[section]
\newtheorem{definition}{Definition}[section]
\newtheorem{conjecture}{Conjecture}[section]

\renewcommand{\theequation}{\thesection.\arabic{equation}}

\newcommand{\BR}{\mathbb R}
\newcommand{\BQ}{\mathbb Q}

\def\charf {\mbox{{\text 1}\kern-.24em {\text l}}}
\newcommand{\BC}{\mathbb C}
\newcommand{\BZ}{\mathbb Z}

\font\Bbb=msbm10 
\begin{document}
\thispagestyle{plain}
 \markboth{}{}
\small{\addtocounter{page}{0} \pagestyle{plain}
\noindent{\scriptsize KYUNGPOOK Math. J. 00(0000), 000-000}\\
\noindent{\scriptsize https://doi.org/10.5666/KMJ.0000.00.0.000}\\
\noindent {\scriptsize pISSN 1225-6951 \,\,\,\,\,\, eISSN 0454-8124}\\
\noindent {\scriptsize $\copyright$ Kyungpook Mathematical Journal}
\vspace{0.2in}\\
\noindent{\large \bf Survey of the Arithmetic and Geometric Approach to the Schottky Problem}
\footnote{{}\\ \\[-0.7cm]
Received March 7, 2023; accepted May 7, 2023.
This work was supported by the Max-Planck-Institut f{\"u}r Mathematik in Bonn.}
\vspace{0.15in}\\
\noindent{\sc Jae-Hyun Yang}
\newline
{\it Department of Mathematics, Inha University, Incheon 22212, Republic of Korea\\
e-mail} : {\verb|jhyang@inha.ac.kr or jhyang8357@gmail.com|}
\vspace{0.15in}\\
{\footnotesize {\sc Abstract.}
In this article, we discuss and survey the recent progress towards the Schottky problem,
and make some comments on the relations between the Andr{\'e}-Oort conjecture, Okounkov convex bodies, Coleman's conjecture, stable modular forms, Siegel-Jacobi spaces, stable Jacobi forms and the Schottky problem.
}
\vspace{0.2in}
\pagestyle{myheadings}
 \markboth{\headd Jae-Hyun Yang$~~~~~~~~~~~~~~~~~~~~~~~~~~~~~~~~~~~~~~~~~~~~~~~~~~~~~~\,$}
 {\headd $~~~~~~~~~~~~~~~~~~~~~~$Arithmetic and Geometric Approach to the Schottky Problem}

\newcommand\tr{\triangleright}
\newcommand\al{\alpha}
\newcommand\be{\beta}
\newcommand\g{\gamma}
\newcommand\gh{\Cal G^J}
\newcommand\G{\Gamma}
\newcommand\de{\delta}
\newcommand\e{\epsilon}
\newcommand\lam{\lambda}
\newcommand\z{\zeta}
\newcommand\vth{\vartheta}
\newcommand\vp{\varphi}
\newcommand\om{\omega}
\newcommand\p{\pi}
\newcommand\la{\lambda}
\newcommand\lb{\lbrace}
\newcommand\lk{\lbrack}
\newcommand\rb{\rbrace}
\newcommand\rk{\rbrack}
\newcommand\s{\sigma}
\newcommand\w{\wedge}
\newcommand\fgj{{\frak g}^J}
\newcommand\lrt{\longrightarrow}
\newcommand\lmt{\longmapsto}
\newcommand\lmk{(\lambda,\mu,\kappa)}
\newcommand\Om{\Omega}
\newcommand\ka{\kappa}
\newcommand\ba{\backslash}
\newcommand\ph{\phi}
\newcommand\M{{\Cal M}}
\newcommand\bA{\bold A}
\newcommand\bH{\bold H}
\newcommand\D{\Delta}

\newcommand\cP{\Cal P}

\newcommand\cH{\Cal H}

\newcommand\pa{\partial}

\newcommand\pis{\pi i \sigma}
\newcommand\sd{\,\,{\vartriangleright}\kern -1.0ex{<}\,}
\newcommand\wt{\widetilde}
\newcommand\fg{\frak g}
\newcommand\fk{\frak k}
\newcommand\fp{\frak p}
\newcommand\fs{\frak s}
\newcommand\fh{\frak h}
\newcommand\Cal{\mathcal}

\newcommand\fn{{\frak n}}
\newcommand\fa{{\frak a}}
\newcommand\fm{{\frak m}}
\newcommand\fq{{\frak q}}
\newcommand\CP{{\mathcal P}_n}
\newcommand\Hnm{{\mathbb H}_n \times {\mathbb C}^{(m,n)}}
\newcommand\BD{\mathbb D}
\newcommand\BH{\mathbb H}
\newcommand\CCF{{\mathcal F}_n}
\newcommand\CM{{\mathcal M}}
\newcommand\Gnm{\Gamma_{n,m}}
\newcommand\Cmn{{\mathbb C}^{(m,n)}}
\newcommand\Yd{{{\partial}\over {\partial Y}}}
\newcommand\Vd{{{\partial}\over {\partial V}}}

\newcommand\Ys{Y^{\ast}}
\newcommand\Vs{V^{\ast}}
\newcommand\LO{L_{\Omega}}
\newcommand\fac{{\frak a}_{\mathbb C}^{\ast}}

\setcounter{section}{1}
\setcounter{theorem}{0}
\noindent{\bf 1. Introduction}
\setcounter{equation}{0}
\renewcommand{\theequation}{1.\arabic{equation}}
\vspace{0.1in}\\
\indent
For a positive integer $g$, we let
\begin{equation*}
\BH_g =\left\{ \tau\in\BC^{(g,g)}\,|\ \tau =\,{}^t\tau, \ {\rm Im}\,\tau>0 \right\}
\end{equation*}
be the Siegel upper half plane of degree $g$
and let
$$Sp(2g,\BR)=\{ M\in \BR^{(2g,2g)}\ \vert \ ^t\!MJ_gM= J_g\}$$
be the symplectic group of degree $g$, where $F^{(k,l)}$ denotes
the set of all $k\times l$ matrices with entries in a commutative
ring $F$ for two positive integers $k$ and $l$, $^t\!M$ denotes
the transposed matrix of a matrix $M$ and
$$J_g=\begin{pmatrix} 0&I_g\\
                   -I_g&0\end{pmatrix}.$$
Then $Sp(2g,\BR)$ acts on $\BH_g$ transitively by
\begin{equation}
M\cdot\tau=(A\tau+B)(C\tau+D)^{-1},
\end{equation} where $M=\begin{pmatrix} A&B\\
C&D\end{pmatrix}\in Sp(2g,\BR)$ and $\Om\in \BH_n.$ Let
$$\G_g=Sp(2g,\BZ)=\left\{ \begin{pmatrix} A&B\\
C&D\end{pmatrix}\in Sp(2g,\BR) \,\big| \ A,B,C,D\
\textrm{integral} \right\}$$ be the Siegel modular group of
degree $g$. This group acts on $\BH_g$ properly discontinuously.

\vskip 0.35cm
Let ${\Cal A}_g:=\G_g\backslash \BH_g$ be the Siegel modular variety of degree $g$, that is, the moduli space of $g$-dimensional principally polarized abelian varieties,
and let ${\Cal M}_g$ be the the moduli space of projective curves of genus $g$. Then according to Torelli's theorem, the Jacobi mapping
\begin{equation}
T_g:{\Cal M}_g \lrt {\Cal A}_g
\end{equation}
defined by
\begin{equation*}
C \longmapsto J(C):= {\rm the\ Jacobian\ of}\ C
\end{equation*}
is injective. The Jacobian locus $J_g:=T_g({\Cal M}_g)$ is a $(3g-3)$-dimensional subvariety of ${\Cal A}_g$

\vskip 0.35cm
The Schottky problem is to characterize the Jacobian locus or its closure ${\bar J}_g$ in ${\Cal A}_g.$ At first this problem had been investigated from
the analytical point of view : to find explicit equations of $J_g$ (or ${\bar J}_g$) in ${\Cal A}_g$ defined by Siegel modular forms on $\BH_g$, for example,
polynomials in the theta constant $\theta \left[ \begin{matrix} \epsilon \\ \delta \end{matrix} \right](\tau,0)$ (see Definition (2.4)) and their derivatives.
The first result in this direction was due to Friedrich Schottky \cite{Sch} who gave the simple and beautiful equation satisfied by the theta constants of Jacobians of dimension 4.
Much later the fact that this equation characterizes the Jacobian locus $J_4$ was proved by J. Igusa \cite{I2} (see also E. Freitag \cite{Fr4} and Harris-Hulek \cite{H-H}). Past decades there has been some progress on the characterization of Jacobians by some mathematicians. Arbarello and De Concini \cite{A-D1} gave a set of such equations defining ${\bar J}_g.$ The Novikov conjecture which states that a theta function satisfying the Kadomtsev-Petviasvili (briefly, K-P) differential equation is the theta function of a Jacobian was proved by T. Shiota \cite{Sh}. Later the proof of the above Novikov conjecture was simplified by Arbarello and De Concini \cite{A-D2}. Bert van Geeman \cite{vG} showed that ${\bar J}_g$ is an irreducible component of the subvariety of ${\Cal A}^{\rm Sat}_g$ defined by certain equations. Here ${\Cal A}^{\rm Sat}_g$ is the Satake compactification of ${\Cal A}_g$. I. Krichever \cite{K2} proved that the
existence of one trisecant line of the associated Kummer variety characterizes Jacobian varieties among principally polarized abelian varieties.

\vskip 0.35cm
S.-T. Yau and Y. Zhang \cite{YZ} obtained the interesting results about asymptotic behaviors of logarithmical canonical line bundles on toroidal compactifications of the Siegel modular varieties. Working on log-concavity of multiplicities in representation theory, A. Okounkov
\cite{OK1, OK2} showed that one could associate a convex body to a linear system on a projective variety, and use convex geometry to study such linear systems. Thereafter R. Lazarsfeld and M. Mustat\u{a} \cite{LM} developed the theory of Okounkov convex bodies associated to linear series systematically. E. Freitag \cite{Fr2} introduced the concept of stable modular forms to investigate the geometry of the Siegel modular varieties. In 2014, using stable modular forms, G. Codogni and N.~I. Shepherd-Barron \cite{Cod-Sh} showed there is no stable Schottky-Siegel forms. We recall that Schottky-Siegel forms
are scalar-valued Siegel modular forms vanishing on the Jacobian locus. Recently G. Codogni \cite{Cod} found the ideal of stable equations of the hyperelliptic locus.
About twenty years ago the author \cite{YJH3, YJH9} introduced the notion of stable Jacobi forms to try to study the geometry of the universal abelian varieties.
In this paper, we discuss the relations among logarithmical line bundles on toroidal compactifications, the Andr{\'e}-Oort conjecture, Okounkov convex bodies, Coleman's conjecture, Siegel-Jacobi spaces, stable Schottky-Siegel forms, stable Schottky-Jacobi forms and the Schottky problem.


\vskip 0.35cm
This article is organized as follows. In Section 2, we briefly survey some known approaches to the Schottky problem and some results so far obtained concerning the characterization of Jacobians.
In Section 3, we briefly describe the results of Yau and Zhang concerning the behaviors of logarithmical canonical line bundles on toroidal compactifications of the Siegel modular varieties.
In Section 4, we review some recent progress on the Andr{\'e}-Oort conjecture.
In Section 5, we review the theory of Okounkov convex bodies associated to linear series (cf.~\cite{CHPW,LM}).
In Section 6, we discuss the relations among logarithmical line bundles on toroidal compactifications, the Andr{\'e}-Oort conjecture, Okounkov convex bodies, Coleman's conjecture and the Schottky problem.
In the final section we give some remarks and propose some open problems about the relations among the Schottky problem, the Andr{\'e}-Oort conjecture, Okounkov convex bodies, stable Schottky-Siegel forms, stable Schottky-Jacobi forms and the geometry of the Siegel-Jacobi space. We define the notion of stable Schottky-Jacobi forms and the concept of stable Jacobi equations for the universal hyperelliptic locus.
In Appendix A, we survey some known results about subvarieties of the Siegel modular variety.
In Appendix B, we review recent results concerning an extension of the Torelli map to a toroidal compactification of the Siegel modular variety.
In Appendix C, we describe why the study of singular modular forms is closely related to that of the geometry of the Siegel modular variety.
In Appendix D, we briefly talk about singular Jacobi forms.  Finally in Appendix E, we review the concept of stable Jacobi forms introduced by the author and relate the study of stable Jacobi forms to that of the geometry of the universal abelian variety.
\vskip 0.3cm
Finally the author would like to mention that he tried to write this article in another new perspective concerning the Schottky problem different from that of other mathematicians. The list of references in this article is by no means complete though we have strived to give as many references as possible. Any inadvertent omissions of references related to the contents in this paper will be the author's fault.

\vskip 0.51cm \noindent {\bf Notations:} \ We denote by
$\BQ,\,\BR$ and $\BC$ the field of rational numbers, the field of
real numbers and the field of complex numbers respectively. We
denote by $\BZ$ and $\BZ^+$ the ring of integers and the set of
all positive integers respectively. $\BR^+$ denotes the set of all positive real numbers.
$\BZ_+$ and $\BR_+$ denote the set of all nonnegative integers and the set of all nonnegative real numbers respectively. The symbol ``:='' means that the expression on the right is the definition of that on the left. For two positive integers $k$ and $l$, $F^{(k,l)}$ denotes the set
of all $k\times l$ matrices with entries in a commutative ring
$F$. For a square matrix $A\in F^{(k,k)}$ of degree $k$,
$\rm{tr}(A)$ denotes the trace of $A$. For any $M\in F^{(k,l)},\
^t\!M$ denotes the transpose of a matrix $M$. $I_n$ denotes the
identity matrix of degree $n$. For $A\in F^{(k,l)}$ and $B\in
F^{(k,k)}$, we set $B[A]=\,^tABA.$ For a complex matrix $A$,
${\overline A}$ denotes the complex {\it conjugate} of $A$. For
$A\in \BC^{(k,l)}$ and $B\in \BC^{(k,k)}$, we use the abbreviation
$B\{ A\}=\,^t{\overline A}BA.$ For a number field $F$, we denote by ${\Bbb A}_{F,f}$ the ring of finite ad{\'e}les of $F.$

\vspace{0.2in}
%
%
%
%
%
%
%
%
%
%
%
%
%
\setcounter{section}{2}
\setcounter{theorem}{0}
\setcounter{equation}{0}
\renewcommand{\theequation}{2.\arabic{equation}}
\noindent{\bf 2. Some Approaches to the Schottky Problem}
\vspace{0.1in}

\newcommand\POB{ {{\partial}\over {\partial{\overline \Omega}}} }
\newcommand\PZB{ {{\partial}\over {\partial{\overline Z}}} }
\newcommand\PX{ {{\partial}\over{\partial X}} }
\newcommand\PY{ {{\partial}\over {\partial Y}} }
\newcommand\PU{ {{\partial}\over{\partial U}} }
\newcommand\PV{ {{\partial}\over{\partial V}} }
\newcommand\PO{ {{\partial}\over{\partial \Omega}} }
\newcommand\PZ{ {{\partial}\over{\partial Z}} }

\indent
Before we survey some approaches to the Schottky problem, we provide some notations and definitions. Most of the materials in this section can be found in \cite{G2}. We refer to \cite{B-D1, DO1, F, G2, Muf5, Ran2} for more details and discussions on the Schottky problem.
\vskip 0.26cm
In this section, we let $g$ be a fixed positive integer. For a positive integer $\ell$, we define the principal level $\ell$ subgroup
$$ \G_g(\ell): = \left\{ \gamma\in Sp(2g,\BZ)\,|\ \gamma \equiv I_{2g}\ ({\rm mod}\,\ell) \,\right\}.$$
and define the theta level $\ell$ subgroup
$$ \G_g(\ell,2\ell):= \left\{  \begin{pmatrix} A & B \\ C & D \end{pmatrix}\in \G_g(\ell)\,\Big| \ {\rm diag}(\,{}^t\!AB)\equiv {\rm diag}(\,{}^tCD) \equiv 0\ ({\rm mod}\,\ell) \,\right\}.$$
We let
\begin{equation*}
{\mathcal A}_g(\ell):=\G_g(\ell)\ba \BH_g\quad {\rm and}\quad {\mathcal A}_g(\ell,2\ell):=\G_g(\ell,2\ell)\ba \BH_g.
\end{equation*}

\begin{definition}
{\rm
(\cite[pp.\,49-50]{I1},~\cite[p.\,123]{Muf3},~ \cite[p.\,862]{YJH7}~or~\cite[p.\,127]{YJH16-1})
Let $\ell$ a positive integer. For any $\epsilon$ and $\delta$ in ${\frac 1\ell}\BZ^g/\BZ^g,$ we define the {\sf theta function with characteristics} $\epsilon$ and $\delta$ by
\begin{equation}
\theta \left[ \begin{matrix} \epsilon \\ \delta \end{matrix} \right](\tau,z):= \sum_{N\in\BZ^g} e^{\pi i\{(N+\epsilon)\,\tau\,{}^t(N+\epsilon)\,+\, 2\,(N+\epsilon)\,{}^t(z+\delta)\}},
\quad (\tau,z)\in \BH_g\times \BC^g.
\end{equation}
The Riemann theta function $\theta(\tau,z)$ is defined to be
\begin{equation}
\theta (\tau,z):=\theta \left[ \begin{matrix} 0 \\ 0 \end{matrix} \right](\tau,z),\quad (\tau,z)\in \BH_g\times \BC^g.
\end{equation}
}
\end{definition}

For each $\tau\in \BH_g,$ we have the transformation behavior
\begin{equation}
\theta (\tau,z+ a\tau +b):=e^{-\pi i (a\tau\,{}^t\!a\, + \,2\,a {}^t\!z)} \theta (\tau,z)\quad {\rm for\ all}\ a,b\in \BZ^g.
\end{equation}
The function
\begin{equation}
\theta \left[ \begin{matrix} \epsilon \\ \delta \end{matrix} \right](\tau) :=\theta \left[ \begin{matrix} \epsilon \\ \delta \end{matrix} \right](\tau,0),\quad \tau\in \BH_g
\end{equation}
is called the {\sf theta constant} of {\it order} $\ell.$ It is known that the theta constants $\theta \left[ \begin{matrix} \epsilon \\ \delta \end{matrix} \right](\tau)$ of order $\ell$
are Siegel modular forms of weight ${\frac 12}$ for $\G_g (\ell, 2\ell)$\,\cite[p.\,200]{Muf3}.

\vskip 0.26cm
For a fixed $\tau\in\BH_g$, we let $\Lambda_\tau:= \BZ^g \tau+\BZ^g$ be the lattice in $\BC^g.$ According to the formula (2.3), the zero locus $\left\{ z\in\BC^g\,|\ \theta (\tau,z)=0\, \right\}$
is invariant under the action of the lattice $\Lambda_\tau$ on $\BC^g$, and thus descends to a well-defined subvariety $\Theta_\tau \subset A_\tau:=\BC^g/ \Lambda_\tau.$
In fact $A_\tau$ is a principally polarized abelian variety with ample divisor $\Theta_\tau.$

\begin{definition}
{\rm
For $\epsilon\in {\frac 12}\BZ^g/\BZ^g$ the theta function of the second order with characteristic $\epsilon$ is defined to be
\begin{equation}
\Theta [\epsilon](\tau,z):= \theta \left[ \begin{matrix} 2\epsilon \\ 0 \end{matrix} \right](2\tau,2z), \quad (\tau,z)\in \BH_g\times \BC^g.
\end{equation}
We define the {\sf theta constant} of the second order to be
\begin{equation}
\Theta [\epsilon](\tau):= \theta \left[ \begin{matrix} 2\epsilon \\ 0 \end{matrix} \right](2\tau,0)=\Theta [\epsilon](\tau,0),
\quad \tau\in \BH_g.
\end{equation}
}
\end{definition}
Then we see that $\Theta [\epsilon](\tau)$ is a Siegel modular form of weight ${\frac 12}$ for $\G_g(2,4).$

\vskip 0.35cm
We have the following results.
\begin{theorem}
{\rm (Riemann's bilinear addition formula)\,\cite[p.\,139]{I1}}\\ \noindent
\begin{equation}
\left( \theta \left[ \begin{matrix} \epsilon \\ \delta \end{matrix} \right](\tau,z) \right)^2= \sum_{\sigma\in {\frac 12}\BZ^g/\BZ^g} (-1)^{4\,{}^t\!\sigma\delta} \Theta[\sigma+\epsilon](\tau,0)\cdot \Theta[\sigma](\tau,z).
\end{equation}
\end{theorem}

\begin{theorem} For $\ell\geq 2,$ the map
\begin{equation}
\Phi_{\ell}: {\mathcal A}_g(2\ell, 4\ell) \lrt {\mathbb P}^N(\BC),\quad N:=\ell^{2g}-1
\end{equation}
defined by
\begin{equation*}
\Phi_{\ell}(\tau):= \left\{ \theta \left[ \begin{matrix} \epsilon \\ \delta \end{matrix} \right](\tau)\,\Big|\ \epsilon,\delta
\in {\frac 1\ell}\BZ^g/\BZ^g\, \right\}
\end{equation*}
is an embedding.
\end{theorem}
\vskip 0.23cm\noindent
{\it Proof.} See Igusa \cite{I1} for $\ell=4n^2$, and Salvati Manni \cite{S0} for $\ell\geq 2.$
 \hfill $\Box$

\begin{remark}
{\rm
We consider the theta map
\begin{equation}
Th: {\mathcal A}_g(2, 4) \lrt {\mathbb P}^{2^g -1}(\BC)
\end{equation}
defined by
\begin{equation*}
Th (\tau):=\left\{ \Theta[\epsilon](\tau)\,\big|\ \epsilon\in {\frac 12}\BZ^g/\BZ^g\, \right\}.
\end{equation*}
We observe that according to Theorem 2.1, $\Phi_2(\tau)$ can be recovered uniquely up to signs from $\Theta(\tau).$
Since $\Phi_2$ is injective on ${\mathcal A}_g(4,8),$ the theta map $Th$ is finite-to-one on ${\mathcal A}_g(2,4).$
In fact, it is known that the theta map $Th$ is generically injective, and it is conjectured that $Th$ is an
embedding.
}
\end{remark}

Now we briefly survey some approaches to the Schottky problem. As mentioned before, most of the following materials in this section comes from a good survey paper \cite{G2}.

\vspace{0.2in}
\noindent
{\bf (A) Classical Approach}

\vspace{0.1in}
For $\tau\in \BH_g$ and a positive integer $\ell\in\BZ^g$,
$$ A_\tau[\ell]:=\left\{ \tau\epsilon +\delta \in A_\tau \,|\ \epsilon,\delta\in {\frac 1\ell}\BZ^g/\BZ^g\right\}$$
denotes the subgroup of $A_\tau$ consisting of torsion points of order $\ell$. For $m=\tau\epsilon +\delta \in A_\tau[\ell]$, we
briefly write
\begin{equation}
\theta_m (\tau,z):=\theta \left[ \begin{matrix} \epsilon \\ \delta \end{matrix} \right](\tau,z).
\end{equation}
\indent
We define the {\sf Igusa\ modular\ form} to be
\begin{equation}
F_g(\tau):= 2^g \sum_{m\in A_\tau [2]} \theta_m^{16}(\tau)-\left( \sum_{m\in A_\tau [2]} \theta_m^8(\tau)\right)^2, \quad \tau\in \BH_g.
\end{equation}
It was proved that $F_g(\tau)$ is a Siegel modular form of weight $8$ for the Siegel modular group $\G_g$ such that when rewritten in terms of
theta constants of the second order using Theorem 2.1,
\medskip
\par
($F_g$--1)\ \ $F_g\equiv 0$\ \ for $g=1,2$\,;
\par
($F_g$--2)\ \ $F_3$ is the defining equation for $\overline{Th(J_3(2,4))}=\overline{Th({\mathcal A}_3(2,4))}\subset {\mathbb P}^7(\BC)$\,;
\par
($F_g$--3)\ \ $F_4$ is the defining equation for $\overline{Th(J_4(2,4))}\subset \overline{Th({\mathcal A}_4(2,4))}\subset {\mathbb P}^{15}(\BC)$.

\vskip 0.35cm
For more detail, we refer to \cite{Fr4, I2, Sch} for the case $g=4$ and refer to \cite{Ran2} for the case $g=5.$
For $g\geq 5$, no similar solution is known or has been proposed.

\begin{theorem}
If $g\geq 5,$ then $F_g$ does not vanish identically on $J_g.$ In fact, the zero locus of $F_5$ on $J_5$ is the locus of trigonal curves.
\end{theorem}
The above theorem was proved by Grushevsky and Salvati Manni \cite{G-SM3}.

\vspace{0.2in}
\noindent
{\bf(B) The Schottky-Jung Approach}

\begin{definition}
{\rm
For an {\'e}tale connected double cover ${\tilde C}\lrt C$ of a curve $C\in {\mathcal M}_g$ (such a curve is given by a two-torsion point $\eta(\not =0)\in J(C)[2]$)
we define the {\sf Prym variety} to
$$ Prym({\tilde C}\lrt C):=Prym(C,\eta):={\rm Ker}_0 (J({\tilde C})\lrt J(C))\in {\mathcal A}_{g-1},$$
where ${\rm Ker}_0$ denotes the connected component of 0 in the kernel and the map $J({\tilde C})\lrt J(C)$ is the norm map corresponding to the cover
${\tilde C}\lrt C$. We denote by ${\mathcal P}_{g-1}\subset {\mathcal A}_{g-1}$ the locus of Pryms of all {\'e}tale double covers of curves in ${\mathcal M}_g.$
The problem of describing ${\mathcal P}_{g-1}$ is called the {\sf Prym-Schottky problem}.
}
\end{definition}

\begin{remark}
{\rm
The restriction of the principal polarization $\Theta_{J({\tilde C})}$ to the Prym gives twice the principal polarization. However this polarization admits a canonical square root, which thus gives a natural principal polarization on the Prym.
}
\end{remark}

\begin{theorem}$(${\sf Schottky-Jung proportionality}$)$ Let $\tau$ be the period matrix of a curve $C$ of genus $g$ and let $\tau_*$ be the period matrix of the Prym for
$\left[ \begin{matrix} 0 & 0 &  \cdots & 0 \\ 1 & 0 & \cdots & 0 \end{matrix}\right]$. Then for any $\epsilon,\delta\in {\frac 12}\BZ^{g-1}/\BZ^{g-1}$
the theta constants of $J(C)$ and of the Prym are related by
\begin{equation}
\left( \theta \left[ \begin{matrix} \epsilon \\ \delta \end{matrix}\right](\tau_*)\right)^2= C\,
\theta \left[ \begin{matrix} 0 & \epsilon \\ 0 & \delta \end{matrix}\right](\tau)\cdot \theta \left[ \begin{matrix} 0 & \epsilon \\1 & \delta \end{matrix}\right](\tau).
\end{equation}
Here the constant $C$ is independent of $\epsilon,\delta.$
\end{theorem}
\noindent
{\it Proof.} See Schottky-Jung \cite{ScJ} and also Farkas \cite{F} for a rigorous proof. \hfill $\Box$

\vspace{0.05in}

\begin{definition}
{\rm
(The Schottky-Jung locus \cite{G2}).
Let $I_{g-1}$ be the defining ideal for the image $\overline{Th({\mathcal A}_{g-1}(2,4))}
\subset {\mathbb P}^{2^{g-1}-1}$ (see Remark 2.1). For any equation $F\in I_{g-1}$, we let
$F_\eta$ be the polynomial equation on $\mathbb P^{2^g-1}$ obtained by using the Schottky-Jung
proportionality to substitute an appropriate polynomial of degree $2$ in terms of theta constants
of $\tau$ for the square of any theta constant of $\tau_*$. Let $S_g^\eta$ be the ideal obtained from $I_{g-1}$ in this way. We define the {\sf big Schottky-Jung locus}
$\mathcal S_g^\eta(2,4)\subset \mathcal A_g(2,4)$ to be the zero locus of $S_g^\eta$. It is not
known that $I_g\subset S_g^\eta$ and thus we define $\mathcal S_g^\eta(2,4)$ within
$\mathcal A_g(2,4)$, and not as a subvariety of the projective space $\mathbb P^{2^g-1}$. We now
define the {\sf small Schottky-Jung locus} to be
\begin{equation}
\mathcal S_g(2,4):=\bigcap_\eta \mathcal S_g^\eta(2,4),
\end{equation}
where $\eta$ runs over the set ${\frac 12}\BZ^{2g}/\BZ^{2g}\backslash \{0\}.$
We note that the action of $\G_g$ permutes the different $\eta$ and the ideals $S_g^\eta$.
Therefore the ideal defining $\mathcal S_g(2,4)$ is $\G_g$-invariant, and the locus
$\mathcal S_g(2,4)$ is a preimage of some $\mathcal S_g\subset \mathcal A_g$ under the
level cover.
}
\end{definition}

\begin{theorem}
{\rm (a)} The Jacobian locus $J_g$ is an irreducible component of the small Schottky-Jung locus ${\mathcal S}_g.$
\par\noindent
{\rm (b)} $J_g(2,4)$ is an irreducible component of the big Schottky-Jung locus ${\mathcal S}_g^{\eta}(2,4)$ for any $\eta.$
\end{theorem}

\noindent
{\it Proof.} The statement (a) was proved by van Geeman \cite{vG} and the statement (b) was proved by Donagi \cite{DO1}.  \hfill $\Box$

\vskip 0.35cm
Donagi \cite{DO2} conjectured the following.
\begin{conjecture}
The small Schottky-Jung locus is equal to the Jacobian locus, that is, ${\mathcal S}_g=J_g.$
\end{conjecture}

\vspace{0.1in}
\noindent
{\bf (C) The Andreotti-Mayer Approach}

\vspace{0.1in}
We let ${\rm Sing}\, \Theta$ be the singularity set of the theta divisor $\Theta$ for a principally polarized abelian variety $(A, \Theta)$.

\begin{theorem}
For a non-hyperelliptic curve $C$ of genus $g$, $\dim ({\rm Sing}\,\Theta_{J(C)})=g-4$, and for a hyperelliptic curve $C$,
$\dim ({\rm Sing}\,\Theta_{J(C)})=g-3$. For a generic principally polarized abelian variety, the theta divisor is smooth.
\end{theorem}
\noindent
{\it Proof.} The proof was given by Andreotti and Mayer \cite{A-M}.
\hfill $\Box$

\begin{definition}
{\rm
We define the {\sf $k$-th Andreotti-Mayler locus} to be
$$ N_{k,\,g}:= \{ (
A,\Theta)\in {\mathcal A}_g\, |\ \dim {\rm Sing}\,\Theta \geq k\}.$$
}
\end{definition}

\begin{theorem}
$N_{g-2,g}={\mathcal A}^{\rm dec}_g.$ Here
$${\mathcal A}^{\rm dec}_g:=\left( \bigcup_{k=1}^{g-1}{\mathcal A}_k\times {\mathcal A}_{g-k}
\right) \subset {\mathcal A}_g$$
denotes the locus of decomposable ppavs $($product of lower-dimensional ppavs$)$ of dimension $g$.
\end{theorem}

\noindent
{\it Proof.} The proof was given by Ein and Lazasfeld \cite{Ei-L}.
\hfill $\Box$

\begin{theorem}
$J_g$ is an irreducible component of $N_{g-4,g}$, and the locus of hyperelliptic Jacobians ${\rm Hyp}_g$ is
an irreducible component of $N_{g-3,g}$.
\end{theorem}

\noindent
{\it Proof.} The proof was given by Andreotti and Mayer \cite{A-M}.
\hfill $\Box$

\begin{theorem}
The Prym locus ${\mathcal P}_g$ is an irreducible component of $N_{g-6,g}.$
\end{theorem}

\noindent
{\it Proof.} The proof was given by Debarre \cite{D}.
\hfill $\Box$

\begin{theorem}
The locus of Jacobians of curves of genus $4$ with a vanishing theta-null is equal to the locus of $4$-dimensional
principally polarized abelian varieties for which the double point singularity of the theta divisor is not ordinary
$($i.e., the tangent cone does not have maximal rank$)$.
\end{theorem}

\noindent
{\it Proof.} See Grushevsky-Salvati Manni \cite{G-SM2} and Smith-Varley \cite{S-V}.
\hfill $\Box$
\vskip 0.35cm\noindent
{\bf Problem.} Can it happen that $N_{k,g}=N_{k+1,g}$ for some $k,g$?

\vspace{0.2in}
\noindent
{\bf (D) The Approach via the K-P Equation}

\vspace{0.1in}
In his study of solutins of nonlinear equations of Korteveg de Vrie type, I. Krichever \cite{K1} proved the following fact\,:
\begin{theorem}
Let $\tau$ be the period matrix of a curve $C$ of genus $g$ and let $\theta(z)$ {\rm (cf.\,(2.2))} be the Riemann theta function of
the Jacobian $J(C)$. Then there exist three vectors $W_1,W_2,W_3$ in $\BC^g$ with $W_1\not=0$ such that, for every $z\in\BC^g$,
the function
\begin{equation}
u(x,y,z;t):= {{\partial^2\ } \over {\partial x^2}} \log \theta (xW_1+yW_2+tW_3+z)
\end{equation}
satisfies the so-called {\sf Kadomstev-Petriashvili\ equation} $($briefly the K-P equation$)$
\begin{equation}
3\,u_{yy}=\left( u_t-3uu_x-2u_{xxx}\right)_x.
\end{equation}
\end{theorem}

S.\,P. Novikov conjectured that $\tau\in\BH_g$ is the period matrix of a curve if and only if the Riemann theta function corresponding to $\tau\in\BH_g$ satisfies the K-P equation in the sense we just explained in Theorem 2.11. Shiota \cite{Sh} proved that the Novikov conjecture is true, following the
work of Mulase \cite{Mu} and Mumford \cite{Muf1}.
Arbarello and De Concini \cite{A-D2} gave another proof of the Novikov conjecture.

\vspace{0.2in}
\noindent
{\bf (E) The Approach via Geometry of the Kummer Variety}

\begin{definition}
{\rm
The map is the embedding given by
\begin{equation}
{\rm Kum}: A_\tau/{\pm 1} \lrt {\mathbb P}^{2^g-1}(\BC),\qquad {\rm Kum}(z)=\left\{ \Theta[\epsilon](\tau,z)\big|\ \epsilon\in {\frac 12}
\BZ^g/\BZ^g \right\}.
\end{equation}
We call the image of ${\rm Kum}$ the {\sf Kummer variety}. Note that the involution $\pm 1$ has $2^{2g}$ fixed points on $A_\tau$ which are
precisely $A_\tau[2]$, and thus the Kummer variety singular at their images in ${\mathbb P}^{2^g-1}(\BC).$
}
\end{definition}

\begin{theorem}
For any points $p_1,p_2,p_3$ of a curve of genus $g$, the following three points on the Kummer variety are collinear\,:
\begin{equation}
{\rm Kum}(p+p_1-p_2-p_3),\ {\rm Kum}(p+p_2-p_1-p_3),\ {\rm Kum}(p+p_3-p_1-p_2).
\end{equation}
\end{theorem}

\noindent
{\it Proof.} See Fay \cite{Fa} and Gunning \cite{Gu1}.
\hfill $\Box$

\begin{theorem}
For any curve $C\in {\mathcal M}_g,$ for any $1\leq k\leq g$ and for any $p_1,\cdots,p_{k+2},q_1,\cdots,q_k\in C$ the $k+2$
points of the Kummer variety
\begin{equation}
{\rm Kum}\left( 2p_j+\sum_{i=1}^kq_i -\sum_{i=1}^{k=2}p_i\right),\quad j=1,\cdots,k+2
\end{equation}
are linearly dependent.
\end{theorem}

\noindent
{\it Proof.} See Gunning \cite{Gu2}.
\hfill $\Box$

\vskip 0.35cm
I. Krichever \cite{K2} gave a complete proof of a conjecture of Welters concerning a condition for an indecomposable principally
polarized abelian variety to be the Jacobian of a curve\,:

\begin{theorem}
Let ${\mathcal A}_g^{\rm ind}:={\mathcal A}_g\backslash {\mathcal A}_g^{\rm dec}$ be the locus of
indecomposable ppavs of dimension $g$.
For a ppav $A\in {\mathcal A}_g^{\rm ind}$, if ${\rm Kum}(A)\subset {\mathbb P}^{2^g-1}$ has one of the following
\vskip 2mm
\indent {\rm(W1)}\ a trisecant line\par
\indent {\rm (W2)}\ a line tangent to it at one point, and intersecting it another point\par
\indent\indent\indent $($this is a semi-degenerate trisecant, when two points of secancy coincides$)$
\par
\indent {\rm (W3)}\ a flex line $($this is a most degenerate trisecant when all three points of
\par
\indent\indent\indent\ \ secancy coincide$)$
\vskip 2mm\noindent
such that none of the points of intersection of this line with the Kummer variety are $A[2]$ $($where ${\rm Kum}(A)$ is singular$)$,
then $A\in J_g.$
\end{theorem}

For the Prym-Schottky problem, it will be natural whether the Kummer varieties of Pryms have any special geometric properties.
Indeed, Beauville-Debarre \cite{B-D2} and Fay \cite{Fa1} obtained the following.
\begin{theorem}
Let $C\in {\mathcal M}_g$. For any $p,p_1,p_2,p_3\in {\tilde C}\lrt Prym({\tilde C}\lrt C)$ on the Abel-Prym curve  the following
four points of the Kummer variety
\begin{eqnarray*}
{\rm Kum}(p+p_1+p_2+p_3),&\quad {\rm Kum}(p+p_1-p_2-p_3),\\
{\rm Kum}(p+p_2-p_1-p_3),&\quad {\rm Kum}(p+p_3-p_1-p_2)
\end{eqnarray*}
lie on a $2$-plane in ${\mathbb P}^{2^g-1}(\BC).$
\end{theorem}

A suitable analog of the trisecant conjecture was found for Pryms using ideas of integrable systems by Grushevsky and Krichever
\cite{G-K}. They proved the following.
\begin{theorem}
If for some $A\in {\mathcal A}_g^{\rm ind}$ and some $p,p_1,p_2,p_3\in A$ the quadrisecant condition in {\rm Theorem 2.15} holds, and
moreover there exists another quadrisecant given by {\rm Theorem 2.15} with $p$ replaced by $-p$, then $A\in \tilde{{\mathcal P}_g}$.
\end{theorem}

\vspace{0.1in}


\noindent
{\bf (F) The Approach via the $\G_{00}$ Conjecture}

\begin{definition}
{\rm
Let $(A,\Theta)\in {\mathcal A}_g.$ The {\sf linear system} $\G_{00}\subset |2\Theta|$ is defined to consist of those sections that vanish to order at least 4
at the origin\,:
\begin{equation}
\G_{00}:=\{ f\in H^0(A, 2\Theta)\,|\ {\rm mult}_0 f \geq 4\}.
\end{equation}
We define the base locus
$$ F_A:=\{ x\in A\,|\ s(x)=0\quad {\rm for\ all}\ s\in \G_{00}\}.$$
}
\end{definition}

\begin{theorem}
For any $g\geq 5$ and any $C\in {\mathcal M}_g$, we have on the Jacobian $J(C)$ of $C$ the equality
\begin{equation*}
F_{J(C)}= C-C  =\{ x-y\in J(C)\,|\ x,y\in C\}.
\end{equation*}
\end{theorem}

\noindent
{\it Proof.} The above theorem was proved by Welters \cite{Wel3} set theoretically and also by Izadi scheme-theoretically. Originally Theorem 2.17 was
conjectured by van Geeman and van der Geer \cite{vG-vdG}.
\hfill $\Box$

\vskip 0.35cm
van Geeman and van der Geer \cite{vG-vdG} conjectured the following.
\vskip 0.35cm\noindent
{\bf $\G_{00}$ \ Conjecture.} Let $(A,\Theta)\in {\mathcal A}_g^{\rm ind}.$ If $F_A\neq 0$, then $A\in J_g.$

\begin{definition}
{\rm
Let $(A,\Theta)\in {\mathcal A}_g$. For any curve $\G$ on $A$ and any point $x\in A,$
we define
\begin{equation*}
\varepsilon(A,\G,x):= {{\Theta_. \G}\over {{\rm mult}_x \G} },\qquad \varepsilon(A,x):=\inf_{\G\ni x} \varepsilon (A,\G,x).
\end{equation*}
We define the {\sf Seshadri constant} of $(A,\Theta)$ by
\begin{equation*}
\varepsilon(A):=\varepsilon(A,\Theta):= \inf_{x\in A}\varepsilon(A,x).
\end{equation*}
}
\end{definition}

\begin{theorem}
If the $\G_{00}$ conjecture holds, hyperelliptic Jacobians are characterized by the value of their Seshadri constants.
\end{theorem}

\noindent
{\it Proof.} See O. Debarre \cite{D2}.
\hfill $\Box$

\begin{theorem}
If some $A\in {\mathcal A}_g^{\rm ind},$ the linear dependence
\begin{equation}
\Theta [\epsilon](\tau,z) = c\, \Theta [\epsilon](\tau,0) + \sum_{1\leq a\leq b\leq g} c_{ab}
{{\partial \Theta [\epsilon](\tau,0)}\over {\partial \tau_{ab}} }
\end{equation}
for some $c, c_{ab}\in \BC\ (1\leq a\leq b\leq g)$ and for all $\epsilon\in {\frac 12}\BZ^g/\BZ^g$
holds with ${\rm rank}(c_{ab})=1$, then $A\in J_g.$
\end{theorem}

\noindent
{\it Proof.} See S. Grushevsky \cite{G1}.
\hfill $\Box$

\vspace{0.2in}

\noindent
{\bf (G) Subvarieties of a ppav: Minimal Cohomology Classes}

\vspace{0.1in}
The existence of some {\it special} subvarieties of a ppav $(A,\Theta)\in {\mathcal A}_g$ gives a criterion that $A$ is the
Jacobian of a curve. We start by observing that for the Jacobian $J(C)$ of a curve $C\in {\mathcal M}_g$ we can map the symmetric product
${\rm Sym}^d(C)\,(1\leq d <g)$ to $J(C)={\rm Pic}^{g-1}(C)$ by fixing a divisor $D\in {\rm Pic}^{g-1-d}(C)$ and mapping
\begin{equation}
\Phi_{(d)}:{\rm Sym}^d(C)\lrt J(C),\qquad (p_1,\cdots,p_d)\mapsto D+\sum_{i=1}^d p_i.
\end{equation}
The image $W^d(C)$ of the map $\Phi_{(d)}$ is independent of $D$ up to translation, and we can compute its cohomology class
\begin{equation*}
\left[ W^d(C)\right]={{[\Theta]^d}\over {(g-d)!}}\in H^{2g-2d}(J(C)),
\end{equation*}
where $[\Theta]$ is the cohomology class of the polarization of $J(C)$. One can show that the cohomology class is indivisible in
cohomology with $\BZ$-coefficients, and we thus call this class {\sf minimal}. We note that $W^1(C)\simeq C.$ These subvarieties
$W^d(C)\ (1\leq d<g)$ are very {\it special}.
\vskip 0.25cm
We have the following criterion.
\begin{theorem}
A ppav $(A,\Theta)\in {\mathcal A}_g$ is a Jacobian if and only if there exists a curve $C\subset A$ with
$[C]={{[\Theta]^{g-1}}\over {(g-1)!}}\in H^{2g-2}(J(C))$ in which case $(A,\Theta)=J(C).$
\end{theorem}

\noindent
{\it Proof.} See Matsusaka \cite{Ma} and Ran \cite{Ran1}.
\hfill $\Box$

\vskip 0.25cm
Debarre \cite{D1} proved that $J_g$ is an irreducible component of the locus of ppavs for which there is a subvariety of the minimal cohomology class. He conjectured the following.
\begin{conjecture}
If a ppav $(A,\Theta)\in {\mathcal A}_g$ has a $d$-dimensional subvariety of minimal class, then it is either the Jacobian of a curve or a five dimensional intermediate Jacobian of a cubic threefold.
\end{conjecture}

This approach to the Schottky problem gives a complete geometric solution to the weaker version of the problem\,: determining whether a given ppav is the Jacobian of a given curve.

\vspace{0.2in}
%
%
%
%
%
%
%
%
%
%
%
%
%
\setcounter{section}{3}
\setcounter{theorem}{0}
\setcounter{equation}{0}
\renewcommand{\theequation}{3.\arabic{equation}}
\noindent{\bf 3. Logarithmical Canonical Line Bundles on Toroidal Compactifications of the Siegel Modular Varieties}
\vspace{0.1in}\\
\indent
In this section, we review the interesting results obtained by S.-T. Yau and Y. Zhang \cite{YZ} concerning the asymptotic behaviors of the logarithmical canonical line bundle on
a toroidal compactification of the Siegel modular variety.

\vskip 0.35cm
Let $\Gamma$ be a neat arithmetic subgroup of $\G_g.$ Let $\mathcal{A}_{g,\G}:=\G\backslash \BH_g$ and $\overline{\mathcal A}_{g,\G}$ be the toroidal compactification of
$\mathcal{A}_{g,\G}$ constructed by a $GL(g,\BZ)$-admissible family of polyhedral decompositions $\Sigma_{{\mathcal F}_0}$ of the cones. Here ${\mathcal F}_0$ denotes
the standard minimal cusps of $\BH_g$. $\overline{\mathcal A}_{g,\G}$ is an algebraic space, but a projective variety in general. Y.-S. Tai proved that if $\Sigma_{tor}^\G$
is projective (see \cite[Chapter IV, Corollary 2.3, p.\,200]{AMRT}), then $\overline{\mathcal A}_{g,\G}$ is a projective variety. It is known that $\overline{\mathcal A}_{g,\G}$
is the unique Hausdorff analytic variety containing $\mathcal{A}_{g,\G}$ as an open dense subset (cf. \cite{AMRT}).

\vskip 0.21cm
Assume the boundary divisor $D_{\infty,\G}$:=$\overline{\mathcal A}_{g,\G}\backslash 
\mathcal{A}_{g,\G}$ is simple normal crossing.
We put $N=g(g+1)/2.$ For each irreducible component $D_i$ of 
$D_{\infty,\G}=\bigcup_{j}D_j$, let $s_i$ a global section of 
the line bundle $[D_i]$ defining $D_i$. Let $\sigma_{\rm max}$ be an arbitrary top-dimensional cone in $\Sigma_{{\mathcal F}_0}$ and renumber all components $D_i's$ of $D_{\infty,\G}$ such that $D_1,
\cdots,D_N$ correspond to the edges of $\sigma_{\rm max}$ with marking order. Yau and Zhang 
\cite[Theorem 3.2]{YZ} showed that the volume form $\Phi_{g,\G}$ on $\mathcal{A}_{g,\G}$ may be written by
\begin{equation*}
  \Phi_{g,\G}={ { 2^{N-g}\,{\rm Vol}_\G (\sigma_{\rm max})^2\,
  d{\mathcal V}_g }\over { \left( \prod_{j=1}^{N} \|s_j\|^2\right)
  F_{\sigma_{\rm max}}^{g+1}(\log \| s_1 \|_1,\cdots,
  \log \| s_N \|_N) } },
\end{equation*}
where $d{\mathcal V}_g$ is a continuous volume form on a partial compactification $\mathcal U_{\sigma_{\rm max}}$ of $\mathcal{A}_{g,\G}$ with
$\mathcal{A}_{g,\G}\subset \mathcal U_{\sigma_{\rm max}} \subset
\overline{\mathcal A}_{g,\G}$, and each $\| \cdot \|_j$ is a suitable  Hermitian metric of the line bundle $[D_j]$ on
$\overline{\mathcal A}_{g,\G}\ (1\leq j \leq N)$ and
$F_{\sigma_{\rm max}}\in\BZ[x_1,\cdots,x_N]$ is a homogeneous polynomial of degree $g$. Moreover the coefficients of
$F_{\sigma_{\rm max}}$ depends only on both $\G$ and
$\sigma_{\rm max}$ with marking order of edges.
Using the above volume form formula they showed that the unique invariant K{\"a}hler-Einstein metric on $\mathcal{A}_{g,\G}$ endows some restraint combinatorial conditions for all smooth toroidal compactifications of $\mathcal{A}_{g,\G}$.

\vskip 2mm
Let $E_1,\cdots,E_d$ be any different irreducible
components of the boundary divisor $D_{\infty,\G}$ such that $\bigcap_{k=1}^d E_k\neq\emptyset$. Let $K_{g,\G}$ be the canonical line bundle on $\overline{\mathcal A}_{g,\G}$.
Yau and Zhang \cite{YZ} also proved the following facts (a) and (b):

\vskip 0.212cm (a) Let $i_1,\cdots,i_d\in \BZ^+.$ If $d \geq g-1$ and $N-\sum_{k=1}^d i_k >2$ (or if $d\geq g$ and $N-\sum_{k=1}^d i_k =1)$,\\
\indent\indent\ \ then we have
\begin{equation*}
\left( K_{g,\G}+D_{\infty,\G} \right)^{N-\sum_{k=1}^d i_k}
\cdot E_1^{i_1}\cdots E_d^{i_d}=0.
\end{equation*}

\vskip 0.212cm (b) $K_{g,\G}+D_{\infty,\G}$ is  not {\it ample} on $\overline{\mathcal A}_{g,\G}$.

\vskip 0.212cm They also showed that if $d < g-1,$ then the intersection number
\begin{equation*}
\left( K_{g,\G}+D_{\infty,\G} \right)^{N-d}\cdot E_1\cdots E_d
\end{equation*}
can be expressed explicitly using the above volume form formula. The proofs of (a) and (b) can be found in \cite[Theorem 4.15]{YZ}.

\vspace{0.2in}

\setcounter{section}{4}
\setcounter{theorem}{0}
\setcounter{equation}{0}
\renewcommand{\theequation}{4.\arabic{equation}}
\noindent{\bf 4. Brief Review on the Andr{\'e}-Oort Conjecture}
\vspace{0.1in}
\indent

In this section we review recent progress on the Andr{\'e}-Oort conjecture quite briefly.
\begin{definition}
{\rm
Let $(G,X)$ be a Shimura datum and let $K$ be a compact open subgroup of $G({\mathbb A}_f).$ We let
\begin{equation*}
Sh_K(G,X):=G(\BQ)\ba X \times G({\mathbb A}_f)/K
\end{equation*}
be the Shimura variety associated to $(G,X).$ An algebraic subvariety $Z$ of the Shimura variety $Sh_K(G,X)$ is said to be
{\sf weakly special} if there exist a Shimura sub-datum $(H, X_H)$ of $(G,X)$, and a decomposition
\begin{equation*}
(H^{\rm ad},X_H^{\rm ad})=(H_1,X_1) \times (H_2,X_2)
\end{equation*}
and $y_2\in X_2$ such that $Z$ is the image of $X_1\times \{ y_2\}$ in $Sh_K(G,X)$. Here $(H^{\rm ad},X_H^{\rm ad})$ denotes
the adjoint Shimura datum associated to $(G,X)$ and $(H_i,X_i)\ (i=1,2)$ are Shimura data. In this definition,
a weakly special subvariety is said to be {\sf special} if it contains a special point and $y_2$ is special.
}
\end{definition}

Andr{\'e} \cite{A} and Oort \cite{O} made conjectures analogous to the Manin-Mumford conjecture where the ambient variety is a Shimura
variety (the latter partially motivated by a conjecture of Coleman \cite{Col}). A combination of these has become known as the
Andr{\'e}-Oort conjecture (briefly the A-O conjecture).
\vskip 0.35cm

\noindent
{\bf A-O Conjecture.} Let $S$ be a Shimura variety and let $\Sigma$ be a set of special points in $S$. Then every irreducible component
of the Zariski closure of $\Sigma$ is a special subvariety.

\begin{definition}
{\rm
\cite{Pi1, Pi2}
A {\sf pre-structure} is a sequence $\Sigma =(\Sigma_n:\,n\geq 1)$ where each $\Sigma_n$ is a
collection of subsets of $\BR^n$. A pre-structure $\Sigma$ is called a {\sf structure} over
the real field if, for all $n,m\geq 1 $ with $m\leq n$, the following conditions are satisfied:
\\
\indent (1) $\Sigma_n$ is a Boolean algebra (under the usual set-theoretic operations);\\
\indent (2) $\Sigma_n$ contains every semi-algebraic subset of $\BR^n$;\\
\indent (3) if $A\in \Sigma_m$ and $B\in\Sigma_n$, then $A\times B\in \Sigma_{m+n}$;\\
\indent (4) if $n\geq m$ and $A\in \Sigma_n$, then $\pi_{n,m}(A)\in \Sigma_m$, where
$\pi_{n,m}:\BR^n\lrt \BR^m$ is a coordinate\\
\indent\ \ \ \ \ projection on the first $m$ coordinates.
\vskip 0.3cm
\noindent If $\Sigma$ is a structure, and, in addition,\\
\vskip 0.05cm
\indent (5) the boundary of every set in $\Sigma_1$ is finite,\\
\vskip 0.05cm
\noindent
then $\Sigma$ is called an {\sf $o$-minimal structure} over the real field.
\vskip 0.2cm
If $\Sigma$ is a structure and $Z\subset \BR^n$,
then we say that $Z$ is {\sf definable} in $\Sigma$ if $Z\in \Sigma_n$.
A function $f:A\lrt B$ is {\sf definable}
in a structure $\Sigma$ if its graph is definable, in which case the domain $A$ of $f$ and
image $f(A)$ are also definable by the definition.
If $A,\cdots,f,\cdots$ are sets or functions, then we denote by $\BR_{A,\cdots,f,\cdots}$
the smallest structure containing $A,\cdots,f,\cdots$. By a {\sf definable family of sets}
we mean a definable subset $Z\subset \BR^n\times\BR^m$ which we view as a family of fibres
$Z_y\subset \BR^n$ as $y$ varies over the projection of $Z$ onto $\BR^m$ which is definable,
along with all the fibres $Z_y$. A family of functions is said to be definable if the family
of their graphs is. A {\sf definable set} usually means a definable set in some o-minimal
structure over the real field.
}
\end{definition}

\begin{remark}
{\rm
The notion of a o-minimal structure grew out of work van den Dries \cite{DR1, DR2} on Tarski's problem concerning the decidability of the real ordered field with the exponential function, and
was studied in the more general context of linearly ordered structures
by Pillay and Steinhorn\,\cite{PS}, to whom the term ``o-minimal" (``order-minimal") is due.
}
\end{remark}

In 2011 Pila gave a unconditional proof of the A-O conjecture for arbitrary products of modular curves using the theory of o-minimality.
\begin{theorem}
Let
\begin{equation*}
X=Y_1\times \cdots \times Y_n\times E_1\times\cdots \times E_m \times {\mathbb G}_m^{\ell},
\end{equation*}
where $n,m,\ell\geq 0,\ Y_i=\G_{(i)}\ba \BH_1 (1\leq i\leq n)$ are modular curves corresponding to congruence subgroups $\G_{(i)}$ of
$SL(2,\BZ)$ and $E_j\,(1\leq j\leq m)$ are elliptic curves defined over ${\overline \BQ}$ and ${\mathbb G}_m$ is the multiplicative group.
Suppose $V$ is a subset of $X$. Then $V$ contains only a finite number of maximal special subvarieties.
\end{theorem}

\noindent
{\it Proof.} See Pila \cite[Theorem 1.1]{Pi1}.
\hfill $\Box$

\vskip 0.25cm
In 2013 Peterzil and Starchenko proved the following theorem using the theory of o-minimality.
\begin{theorem}
The restriction of the uniformizing map $\pi:\BH_g\lrt {\mathcal A}_g$ to the classical fundamental domain for the Siegel modular group
$Sp(2g,\BZ)$ is definable.
\end{theorem}

\noindent
{\it Proof.} See Peterzil and Starchenko \cite{P-S1,P-S2}.
\hfill $\Box$

\vskip 0.25cm
In 2014 Pila and Tsimerman gave a conditional proof of the A-O conjecture for the Siegel modular variety $\mathcal A_g.$

\begin{theorem}
If $g\leq 6,$ then the A-O conjecture holds for $\mathcal A_g.$ If $g\geq 7$, the A-O conjecture holds for $\mathcal A_g$ under the
assumption of the Generalized Riemann Hypothesis $($GRH$)$ for CM fields.
\end{theorem}
\noindent
{\it Proof.} See Pila-Tsimerman \cite{PT1,PT2}.
\hfill $\Box$

\vskip 0.25cm
Quite recently using Galois-theoretic techniques and geometric properties of Hecke correpondences, Klingler and Yafaev proved
the A-O conjecture for a general Shimura variety, and independently using Galois-theoretic and ergodic techniques Ullmo and Yafaev
proved the A-O conjecture for a general Shimura variety, under the assumption of the GRH for CM fields or another suitable assumption.
The explicit statement is given as follows.

\begin{theorem}
Let $(G,X)$ be a Shimura datum and $K$ a compact open subgroup of $G({\mathbb A}_f)$. Let $\Sigma$ be a set of special points in
$Sh_K(G,X)$. We make one of the two following assumptions\,:
\vskip 0.25cm
{\rm (1)} Assume the GRH for CM fields.
\vskip 0.15cm
{\rm (2)} Assume that there exists a faithful representation $G\hookrightarrow GL_n$ such that with respect to this representation, the
Mumford-Tate group $MT(s)$ lie in one $GL_n(\BQ)$-conjugacy class as $s$ ranges through $\Sigma.$ Then every irreducible
component of $\Sigma$ in $Sh_K(G,X)$ is a special subvariety.
\end{theorem}

\noindent
{\it Proof.} See Klingler-Yafaev \cite{K-Y} and Ullmo-Yafaev \cite{U-Y}.
\hfill $\Box$


\begin{remark}
{\rm
We refer to \cite{Pi2} for the theory of o-minimality and the A-O conjecture. We also refer to \cite{G} for
the A-O conjecture for mixed Shimura varieties.
}
\end{remark}

\vspace{0.1in}

\setcounter{section}{5}
\setcounter{theorem}{0}
\setcounter{equation}{0}
\renewcommand{\theequation}{5.\arabic{equation}}
\noindent{\bf 5. Okounkov Bodies Associated to Divisors}
\vspace{0.1in}
\indent

In this section, we briefly review the theory of Okounkov convex bodies associated to pseudoeffective divisors on a smooth projective variety.
For more details of this theory, we refer to \cite{CHPW,LM}.

\vskip 0.21cm
Let $X$ be a smooth projective variety of dimension $d$. We fix an admissible flag $Y_\bullet$ on $X$
$$ Y_\bullet :\ X=Y_0 \supset Y_1\supset Y_2 \supset \cdots \supset Y_{d-1} \supset Y_d=\{x\},$$
where each $Y_k$ is a subvariety of $X$ of codimension $k$ which is nonsingular at $x$. We let $\BZ_+$ denote the set of all non-negative integers.
We first assume that $D$ is a big Cartier divisor on $X$.
For a section $s\in H^0(X,\mathcal O_X (D))\backslash \{0\},$ we define the function
\begin{equation*}
\nu(s)=\nu_{Y_\bullet}(s):= (\nu_1(s),\cdots, \nu_d(s))\in \BZ_+^d
\end{equation*}
as follows:

\vskip 0.21cm
First we set $\nu_1=\nu_1(s):={\rm ord}_{Y_1}(s).$ Using a local equation $f_1$ for $Y_1$ in $X$, we define naturally a section
$$ s_1'=s \otimes f_1^{-\nu_1} \in H^0(X,\mathcal O_X (D-\nu_1 Y_1))$$
that does not vanish along $Y_1$, its restriction $s_1'|_{Y_1}$ defines a nonzero section
$$ s_1:=s_1'|_{Y_1} \in H^0(Y_1,\mathcal O_{Y_1} (D-\nu_1 Y_1)).$$
We now take
$$\nu_2(s):={\rm ord}_{Y_2}(s_1).$$
and continue in this manner to define the remaining $\nu_i(s).$

\vskip 0.21cm
Next we define
$${\rm vect}(|D|)={\rm Im}\left( \nu_{Y_\bullet}:(H^0(X,\mathcal O_X (D))-\{0\})\lrt \BZ^d\right)$$
be the set of valuation vectors of non-zero sections of $\mathcal O_X(D).$ Then we finally set
\begin{equation}
\Delta(D):=\Delta_{Y_\bullet}(D)={\rm closed\ convex\ hull}\left( \bigcup_{m\geq 1}{\frac 1m}\cdot {\rm vect}(|mD|)\right).
\end{equation}
Therefore $\Delta(D)$ is a convex body in $\BR^d$ that is called the {\sf Okounkov body} of $D$ with respect to the fixed flag $Y_\bullet$.
We refer to \cite[\S 1.2]{LM} for some properties and examples of $\Delta(D)$.

\vskip 0.17cm
We recall that a {\sf graded linear series} $W_\bullet (D)=\{ W_m(D)\}_{m\geq 0}$ associated to $D$ consists of subspaces
$$W_m:=W_m(D)\subseteq H^0(X,\mathcal O_X (mD)),\qquad W_0=\BC$$
satisfying the inclusion
$$W_k\cdot W_l\subseteq W_{k+l} \qquad {\rm for\ all}\ k,l\geq 0.$$
Here the product on the left denotes the image of $W_k\otimes W_l$ under the multiplication map
$H^0(X,\mathcal O_X (kD))\otimes H^0(X,\mathcal O_X (lD))\lrt H^0(X,\mathcal O_X ((k+l)D)).$

\begin{definition}
{\rm
(\cite[Definition 1.16] {LM})
Let $W_\bullet$ be a graded linear series on $X$ belonging to a divisor $D$. The {\sf graded\ semigroup} of $W_\bullet$ is defined to be
$$\G (W_\bullet):=\G_{Y_\bullet} (W_\bullet)=\left\{ (\nu_{Y_\bullet}(s),m)\,|\ 0\neq s\in W_m,\ m
\geq 0 \right\}\subseteq \BZ_+^d\times \BZ_+ \subseteq \BZ^{d+1}.$$
Under the above notations, we associate the convex body $\Delta_{Y_\bullet}(W_\bullet)$ of a graded linear series $W_\bullet$ with respect to $Y_\bullet$ on $X$
as follows:
\begin{equation}
\Delta_{Y_\bullet}(W_\bullet):=\sum \left( \G(W_\bullet)\right) \cap \left( \BR_+^d \times \{1\}\right),
\end{equation}
where $\BR_+$ denotes the set of all non-negative real numbers and $\sum \left( \G(W_\bullet)\right)$ denotes the closure of the convex cone in $\BR_+^d\times \BR_+$
spanned by $\G(W_\bullet)$. $\Delta_{Y_\bullet}(W_\bullet)$ is called the {\sf Okounkov body} of $W_\bullet$ with respect to $Y_\bullet$. If $W_\bullet$ is a complete
graded linear series, that is, $W_m=H^0(X,{\mathcal O}(mD))$ for each $m$, then we define
\begin{equation}
\Delta_{Y_\bullet}(D):=\Delta_{Y_\bullet}(W_\bullet).
\end{equation}
}
\end{definition}

\begin{remark}
{\rm
$\Delta_{Y_\bullet}(D)$ depends on the choice of an admissible flag $Y_\bullet$. By the homogeneity of $\Delta_{Y_\bullet}(D)$ (see \cite[Proposition 4.13]{LM}),
we can extend the construction of $\Delta_{Y_\bullet}(D)$ to $\BQ$-divisors $D$ and even to $\BR$-divisors using the continuity of $\Delta_{Y_\bullet}(D)$.
}
\end{remark}

\begin{definition}
{\rm
(\cite[Definition 2.5 and 2.9] {LM})\\
\noindent
(I) We say that a graded linear series $W_\bullet$ satisfies {\sf Condition (B)}
if $W_m\neq 0$ for all\\
\indent  $m\gg 0$, and for all sufficiently large $m$, the rational map
$\phi_m:X--> \mathbb P (W_m)$ \\
\indent defined by $|W_m|$ is birational on its image.  \\
\noindent
(II) We say that a graded linear series $W_\bullet$ satisfies {\sf Condition (C)} if \\
\indent (1) for any $m\gg 0$, there exists an effective divisor $F_m$ such that $A_m:=mD-F_m$ \\
\indent\ \ \ \ \ is ample, and \\
\indent (2) for all sufficiently large $t$, we have
$$ H^0(X, \mathcal O_X(tA_m))\subseteq W_{tm}\subseteq H^0(X, \mathcal O_X(tmA_m)).  $$
If $W_\bullet$ is complete, that is, $W_m=H^0(X, \mathcal O_X(mA_m))$ for all $m\geq 0$ and $D$ is big,\\
\noindent then it satisfies Condition (C).
}
\end{definition}

Lazarsfeld and Mustat\u{a} \cite{LM} proved the following.
\begin{theorem}
Let $X$ be a smooth projective variety of dimension $d$.
Suppose that a graded linear series $W_\bullet$ satisfies Condition {\rm (B)} or Condition {\rm (C)}. Then for any admissible flag $Y_\bullet$ on $X$, we have
\begin{equation*}
\dim\, \Delta_{Y_\bullet}(W_\bullet)=\dim X =d
\end{equation*}
and
\begin{equation*}
{\rm Vol}_{\BR^n}(\Delta_{Y_\bullet}(W_\bullet))= {1\over {d !}} {\rm Vol}_X(W_\bullet),
\end{equation*}
where
\begin{equation*}
 {\rm Vol}_X(W_\bullet):=\lim_{n\lrt \infty}{{\dim W_m}\over {m^d/d!}}.
\end{equation*}
\end{theorem}
\vskip 0.23cm\noindent
{\it Proof.} See \cite[Theorem 2.13]{LM}.
\hfill $\Box$

\vspace{0.05in}

\begin{remark}
{\rm
It is known by Lazarsfeld and Mustat\u{a} (\cite[Proposition 4.1]{LM}) that for a fixed admissible flag $Y_\bullet$ on $X$, if $D$ is big, then $\Delta_{Y_\bullet}(D)$
depends only on the numerical class of $D$. If $D$ is not big, then it is not true (cf.~\cite[Remark 3.13]{CHPW}).
}
\end{remark}

\begin{definition}
{\rm
For a divisor $D$ on $X$, we let
$$ \mathbb{N}(D):=\left\{ m\in \BZ^+\,|\ | \lfloor mD \rfloor |\neq \emptyset \right\}.$$
For $m\in \mathbb{N}(D),$ we let
$$\Phi_{mD}:X ---> \mathbb P^{\dim | \lfloor mD \rfloor |}$$
be the rational map defined by the linear system $| \lfloor mD \rfloor |$. We define the {\sf Iitaka dimension} of $D$ as the following value
$$\kappa (D):=\begin{cases} {\rm max}\{ \dim {\rm Im}(\Phi_{mD})\,|\ m\in \mathbb{N}(D)\} & \quad {\rm if}\  \mathbb{N}(D)\neq \emptyset \\
\indent -\infty & \quad {\rm if}\  \mathbb{N}(D)= \emptyset.\end{cases}$$
}
\end{definition}

\begin{definition}
{\rm
Let $D$ be a divisor on $X$ such that $\kappa(D)\geq 0.$ A subset $U$ of $X$ is called a {\sf Nakayama subvariety} of $D$ if $\kappa(D)=\dim D$ and the natural map
\begin{equation*}
H^0(X,\mathcal O_X(\lfloor mD \rfloor)) \lrt H^0(U,\mathcal O_U(\lfloor mD|_U \rfloor))
\end{equation*}
is injective for every non-negative integer $m$.
}
\end{definition}

\begin{definition}
{\rm
\cite[Definition 3.8]{CHPW}
Let $D$ be a divisor on $X$ such that $\kappa(D)\geq 0.$ The {\sf valuative Okounkov body} $\Delta_{Y_\bullet}^{\rm val}(D)$ associated to $D$ is defined to be
\begin{equation*}
\Delta_{Y_\bullet}^{\rm val}(D):=\Delta_{Y_\bullet}(D)\subset \BR^n,\qquad n=\dim X.
\end{equation*}
For a divisor $D$ with $\kappa(D)=-\infty,$ we define $\Delta_{Y_\bullet}^{\rm val}(D):=\emptyset.$
}
\end{definition}

\begin{remark}
{\rm
If $D$ is big, then $\Delta_{Y_\bullet}^{\rm val}(D)$ coincides with $\Delta_{Y_\bullet}(D)$ for any admissible flag $Y_{\bullet}$ on $X$.
}
\end{remark}

Recently Choi, Hyun, Park and Won \cite{CHPW} showed the following.
\begin{theorem}
Let $D$ be a divisor with $\kappa (D)\geq 0$ on a smooth projective variety $X$ of dimension $n$. Fix an admissible flag $Y_\bullet$ containing a Nakayama subvariety $U$ of $D$
such that $Y_n=\{x\}$ is a general point. Then we have
\begin{equation*}
\dim\, \Delta_{Y_\bullet}^{\rm val}(D)=\kappa (D)
\end{equation*}
and
\begin{equation*}
{\rm Vol}_{\BR^{\kappa(D)}} (\Delta_{Y_\bullet}^{\rm val}(D))= {1\over {\kappa(D) !}} {\rm Vol}_{X|U}(D).
\end{equation*}
\end{theorem}
\vskip 0.23cm\noindent
{\it Proof.} See \cite[Theorem 3.12]{CHPW}.
\hfill $\Box$

\begin{definition}
{\rm
\cite[Definition 3.17]{CHPW}
Let $D$ be a pseudo-effective divisor on a projective variety $X$ of dimension $n$. The {\sf limiting Okounkov body} $\Delta_{Y_\bullet}^{\rm lim}(D)$ of $D$  with respect to
an admissible flag $Y_\bullet$ is defined to be
\begin{equation*}
\Delta_{Y_\bullet}^{\rm lim}(D):=\lim_{\varepsilon\rightarrow 0^+}\Delta_{Y_\bullet}(D+\varepsilon A) \subset \BR^n,
\end{equation*}
where $A$ is any ample divisor on $X$. If $D$ is not a pseudo-effective divisor, then we define $\Delta_{Y_\bullet}^{\rm lim}(D):=\emptyset.$
}
\end{definition}

\begin{definition}
{\rm
\cite[Definition 2.11]{CHPW}
Let $D$ be a divisor on a projective variety $X$ of dimension $d$. We defne the {\sf numerical Iitaka dimension} $\kappa_\nu(D)$
by
\begin{equation*}
\kappa_\nu(D):={\rm max}\left\{  k\in\BZ_+\, \big|\ \limsup_{m\to\infty} {{h^0(X,\mathcal O_X (|mD|+A))}\over {m^k}}>0 \right\}
\end{equation*}
for a fixed ample Cartier divisor $A$ if $h^0(X,\mathcal O_X (|mD|+A))\neq \emptyset$ for infinitely many $m>0$, and we define
$\kappa_\nu(D)=-\infty$ otherwise.
}
\end{definition}

Let $D$ be a pseudo-effective Cartier divisor on a projective variety $X$ of dimension $n$. Let $V\subseteq X$ be a positive volume subvariety of
$D$. Fix an admissible flag $V_\bullet$ on $V$
$$ V_\bullet :\ V=V_0 \supset V_1\supset V_2 \supset \cdots \supset V_{n-1} \supset V_n=\{x\}.$$
let $A$ be an ample Catier divisor on $X$. For each positive integer $k$, we consider the restricted graded linear series
$W_\bullet^k:=W_\bullet(kD+A|V)$ of $kD+A$ along $V$ given by
$$W_m(kD+A|V)=H^0(X|V,m(kD+A))\qquad {\rm for}\ m \geq 0.$$
We define the {\sf restricted limiting Okounkov body} of a Cartier divisor $D$ with respect to a positive volume subvariety $V$ of $D$ as
\begin{equation*}
\Delta_{V_\bullet}^{\rm lim}(D):=\lim_{k\rightarrow \infty}{1\over k}\Delta_{V_\bullet}(W_\bullet^k) \subseteq \BR^{\kappa_\nu(D)}.
\end{equation*}
By the continuity, we can extend this definition for any pseudo-effective $\BR$-divisor.

\begin{definition}
{\rm
Let $D$ be a pseudo-effective divisor on a projective variety $X$ of dimension $n$ with its positive volume subvariety $V\subseteq X$.
We define the {\sf restricted limiting Okounkov body} $\Delta_{V_\bullet}^{\rm lim}(D)$ of $D$  with respect to
an admissible flag $V_\bullet$ to be a closed convex subset
\begin{equation*}
\Delta_{V_\bullet}^{\rm lim}(D):=\lim_{\varepsilon\rightarrow 0^+}\Delta_{V_\bullet}(D+\varepsilon A) \subseteq \BR^{\kappa_\nu(D)}
\hookrightarrow \BR^n,
\end{equation*}
where $A$ is any ample divisor on $X$. If $D$ is not a pseudo-effective divisor, then we define $\Delta_{V_\bullet}^{\rm lim}(D):=\emptyset.$
}
\end{definition}

Recently Choi, Hyun, Park and Won \cite{CHPW} proved the following.
\begin{theorem}
Let $D$ be a pseudo-effective divisor on a projective variety $X$. Fix a positive volume subvariety $V\subseteq X$ of $D$
{\rm (see \cite[Definition 2.13]{CHPW})}. For an
admissible flag $V_\bullet$ of $V$, we have
\begin{equation*}
\dim\, \Delta_{V_\bullet}^{\rm lim}(D)=\kappa_\nu (D)
\end{equation*}
and
\begin{equation*}
{\rm Vol}_{\BR^{\kappa_\nu(D)}} (\Delta_{V_\bullet}^{\rm lim}(D))= {1\over {\kappa_\nu(D) !}} {\rm Vol}_{X|V}^+(D).
\end{equation*}
Here ${\rm Vol}_{X|V}^+(D)$ denotes the augmented restricted volume of $D$ along $V$
{\rm (see \cite[Definition 2.2]{CHPW})} for the
precise definition of ${\rm Vol}_{X|V}^+(D))$.
\end{theorem}

\noindent
{\it Proof.} See \cite[Theorem 3.20]{CHPW}.
\hfill $\Box$

\newpage

\setcounter{section}{6}
\setcounter{theorem}{0}
\setcounter{equation}{0}
\renewcommand{\theequation}{6.\arabic{equation}}
\noindent{\bf 6. The Relations of the Schottky Problem to the Andr{\'e}-Oort Conjecture, Okounkov Bodies and Coleman's Conjecture}
\vspace{0.1in}

In this section, we discuss the relations among logarithmical line bundles on toroidal compactifications, the Andr{\'e}-Oort conjecture, Okounkov convex bodies, Coleman's conjecture and the Schottky problem.

\vskip 0.21cm For $\tau=(\tau_{ij})\in\BH_g,$ we write $\tau=X+i\,Y$
with $X=(x_{ij}),\ Y=(y_{ij})$ real. We put $d\tau=(d\tau_{ij})$ and $d{\overline\tau}=(d{\overline\tau}_{ij})$. We
also put
$$\PO=\,\left(\,
{ {1+\delta_{ij}}\over 2}\, { {\partial}\over {\partial \tau_{ij} }
} \,\right) \qquad\text{and}\qquad \POB=\,\left(\, {
{1+\delta_{ij}}\over 2}\, { {\partial}\over {\partial {\overline
{\tau}}_{ij} } } \,\right).$$ C. L. Siegel \cite{Si1} introduced
the symplectic metric $ds_{g;A}^2$ on $\BH_g$ invariant under the action
(1.1) of $Sp(2g,\BR)$ that is given by
\begin{equation}
ds_{g;A}^2=A\,{\rm tr} (Y^{-1}d\tau\, Y^{-1}d{\overline\tau}),\qquad A\in \BR^+
\end{equation}
and H.
Maass \cite{M2} proved that its Laplacian is given by
\begin{equation}
\Box_{g;A}=\,{4\over A}\,{\rm tr} \left(\,Y\,
{}^{{}^{{}^{{}^\text{\scriptsize $t$}}}}\!\!\!
\left(Y\POB\right)\PO\right).\end{equation}
Here ${\rm tr}(M)$ denotes the trace of a square matrix $M$. And
\begin{equation}
dv_g(\tau)=(\det Y)^{-(g+1)}\prod_{1\leq i\leq j\leq g}dx_{ij}\,
\prod_{1\leq i\leq j\leq g}dy_{ij}\end{equation} is a
$Sp(2g,\BR)$-invariant volume element on
$\BH_g$\,(cf.\,\cite[p.\,130]{Si2}).

\vskip 2mm
Siegel proved the following theorem for the Siegel space $(\BH_g, ds^2_{g;1}).$

\begin{theorem}\,$(${\bf Siegel\,\cite{Si1}}$)$.
{\rm (1)} There exists exactly one geodesic joining two arbitrary points
$\tau_0,\,\tau_1$ in $\BH_g$. Let $R(\tau_0,\tau_1)$ be the
cross-ratio defined by
\begin{equation*}
R(\tau_0,\tau_1)=(\tau_0-\tau_1)(\tau_0-{\overline
\tau}_1)^{-1}(\overline{\tau}_0-\overline{\tau}_1)(\overline{\tau}_0-\tau_1)^{-1}.
\end{equation*}
For brevity, we put $R_*=R(\tau_0,\tau_1).$ Then the symplectic
length $\rho(\tau_0,\tau_1)$ of the geodesic joining $\tau_0$ and
$\tau_1$ is given by
\begin{equation*}
\rho(\tau_0,\tau_1)^2=\s \left( \left( \log { {1+R_*^{\frac 12}
}\over {1-R_*^{\frac 12} } }\right)^2\right),
\end{equation*} where
\begin{equation*}
\left( \log { {1+R_*^{\frac 12} }\over {1-R_*^{\frac 12} }
}\right)^2=\,4\,R_* \left( \sum_{k=0}^{\infty} { {R_*^k}\over
{2k+1}}\right)^2.
\end{equation*}

\noindent {\rm (2)} For $M\in Sp(2g,\BR)$, we set
$${\tilde \tau}_0=M\cdot \tau_0\quad \textrm{and}\quad {\tilde \tau}_1=M\cdot
\tau_1.$$ Then $R(\tau_1,\tau_0)$ and
$R({\tilde\tau}_1,{\tilde\tau}_0)$ have the same eigenvalues.

\vskip 2mm\noindent
\noindent {\rm (3)} All geodesics are symplectic images of the special
geodesics
\begin{equation*}
\alpha(t)=i\,{\rm diag}(a_1^t,a_2^t,\cdots,a_g^t),
\end{equation*}
where $a_1,a_2,\cdots,a_g$ are arbitrary positive real numbers
satisfying the condition
$$\sum_{k=1}^g \left( \log a_k\right)^2=1.$$
\end{theorem}
\noindent The proof of the above theorem can be found in
\cite{M3} or \cite[pp.\,289-293]{Si1}.

\begin{definition}
{\rm
Let $Z$ be an irreducible subvariety of a Shimura variety $Sh_K(G,X)$. Choose a connected component $S$ of $X$
and a class $\eta K\in G({\mathbb A}_f)/K$ such that $Z$ is contained in the image of $S$ in $Sh_K(G,X)$.
We say that $Z$ is a {\sf totally geodesic} subvariety if there is a totally geodesic subvariety $Y \subseteq S$
such that $Z$ is the image of $Y\times \eta K$ in $Sh_K(G,X)$.
}
\end{definition}

B. Moonen \cite{Mo} proved the following fact.
\begin{theorem}
Let $Z$ be an irreducible subvariety of a Shimura variety $Sh_K(G,X)$. Then $Z$ is weakly special if and only if
it is totally geodesic.
\end{theorem}

\noindent
{\it Proof.} See \cite[Theorem 4.3, pp.\,553--554]{Mo}.
\hfill $\Box$

\vspace{0.1in}
In the 1980s Coleman \cite{Col} proposed the following conjecture.
\vskip 0.25cm\noindent
{\bf Coleman's Conjecture.} {\it For a sufficiently large integer $g$, the Jacobian locus $J_g$ contains only a finite number of
special points in ${\mathcal A}_g.$}

\vskip 0.35cm
We also have the following conjecture.

\begin{conjecture}
For a sufficiently large integer $g$, the Jacobian locus $J_g$ cannot contain a non-trivial totally geodesic subvariety.
\end{conjecture}

\begin{remark}
{\rm
Conjecture 6.1 is false for an integer $g\leq 6.$
}
\end{remark}

The stronger version of Conjecture 6.1 is given as follows:
\begin{conjecture}
For a sufficiently large integer $g$, there does not exist a geodesic in ${\mathcal A}_g$ contained in ${\overline J}_g$
and intersecting $J_g.$
\end{conjecture}

\begin{theorem}
Suppose the Andr{\'e}-Oort conjecture and {\rm Conjecture 6.1} hold. Then Coleman's conjecture is true.
\end{theorem}

\noindent
{\it Proof.} Let $g$ be a sufficiently large integer $g$. Suppose $J_g$ contains an infinite set $\Sigma$ of special points.
Then
$$ \Sigma \subset {\overline{\Sigma}}\subset {\overline J}_g \subset {\mathcal A}_g.$$
The truth of the Andr{\'e}-Oort conjecture implies that ${\overline{\Sigma}}$ contains an irreducible special subvariety $Y$.
According to Theorem 6.2, $Y$ is a totally geodesic subvariety of ${\overline J}_g$. From the truth of Conjecture 6.1, we
get a contradiction. Therefore $J_g$ contains only finitely many special points.
\hfill $\Box$

\vspace{0.1in}
Now we propose the following problems.

\vspace{0.1in}
\noindent
{\bf Problem 6.1.} Develop the spectral theory of the Laplace operator $\Box_{g;A,B}$ on $\BH_g$ and $\Gamma\ba \BH_g$ for a congruence subgroup $\G$ of $\G_g$ explicitly.

\vspace{0.1in}
\noindent
{\bf Problem 6.2.} Construct all the geodesics contained in ${\overline J}_g$ with respect to the Siegel's metric $ds^2_{g;A}.$

\vspace{0.1in}
\noindent
{\bf Problem 6.3.} Study variations of $g$-dimensional principally polarized abelian varieties along a geodesic inside $J_g$.

\vspace{0.1in}
\noindent
{\bf Problem 6.4.} Prove the A-O conjecture for ${\mathcal A}_g$ unconditionally.

\vspace{0.1in}
From now on, we will adopt the notations in Section 3.

\vspace{0.1in}
\noindent
{\bf Problem 6.5.} Let $p_{4,\G}:\mathcal A_{4,\G}\lrt \mathcal A_4$ be a covering map and let $J_{4,\G}:=p_{4,\G}^{-1}(J_4).$
Let $\overline{\mathcal A}_{4,\G}$ be a toroidal compactification of $\mathcal A_{4,\G}$ which is projective. Then $J_{4,\G}$ is a
divisor on $\overline{\mathcal A}_{4,\G}$. Compute the Okounkov bodies $\Delta_{Y_\bullet}(J_{4,\G})$, $\Delta_{Y_\bullet}^{\rm val}(J_{4,\G})$ and $\Delta_{Y_\bullet}^{\rm lim}(J_{4,\G})$ explicitly. Describe the relations among $J_4,\ J_{4,\G}$ and these Okounkov bodies explicitly.
Describe the relations between these Okounkov bodies and the $GL(4,\BZ)$-admissible family of polyhedral decompositions defining the toroidal
compactification $\overline{\mathcal A}_{4,\G}$.

\vspace{0.1in}
\noindent
{\bf Problem 6.6.} Assume that a toroidal compactification $\overline{\mathcal A}_{g,\G}$ is a projective variety. Compute the Okoukov convex bodies
$\Delta_{Y_\bullet}(K_{g,\G})$, $\Delta_{Y_\bullet}^{\rm val}(K_{g,\G})$,
$\Delta_{Y_\bullet}^{\rm lim}(K_{g,\G})$, $\Delta_{Y_\bullet}(D_{\infty,\G})$,
$\Delta_{Y_\bullet}^{\rm val}(D_{\infty,\G})$,
$\Delta_{Y_\bullet}^{\rm lim}(D_{\infty,\G})$,
$\Delta_{Y_\bullet}(K_{g,\G}+D_{\infty,\G})$,
$\Delta_{Y_\bullet}^{\rm val}(K_{g,\G}+D_{\infty,\G})$,
$\Delta_{Y_\bullet}^{\rm lim}(K_{g,\G}+D_{\infty,\G})$ explicitly. Describe the relations between these Okounkov bodies and the $GL(g,\BZ)$-admissible family of polyhedral decompositions defining the toroidal
compactification $\overline{\mathcal A}_{g,\G}$.

\vspace{0.1in}
\noindent
{\bf Problem 6.7.} Assume that $g\geq 5.$ Let $p_{g,\G}:\mathcal A_{g,\G}\lrt \mathcal A_g$ be a covering map and let $J_{g,\G}:=p_{g,\G}^{-1}(J_g).$
Assume that $\overline{\mathcal A}_{g,\G}$ is a toroidal compactification of $\mathcal A_{g,\G}$ which is a projective variety.
Let $D_{J,\G}$ be a divisor on $\overline{\mathcal A}_{g,\G}$ containing $J_{g,\G}$. Describe the Okounkov bodies
$\Delta_{Y_\bullet}(D_{J,\G}),\ \Delta_{Y_\bullet}^{\rm val}(D_{J,\G}),\ \Delta_{Y_\bullet}^{\rm lim}(D_{J,\G})$.
Study the relations between $J_{g,\G},\ D_{J,\G}$ and these Okounkov bodies.

\vskip 0.251cm
We have the following diagram:
\begin{equation*}
\begin{tikzcd}[row sep=0.7cm]
J_{g,\Gamma}\ar[d]\ar[r,phantom,"\ensuremath{\subset}"] & {\mathcal A}_{g,\Gamma}
\ar[r,phantom,"\ensuremath{\subset}"]
\ar[d]\ar[d,"p_{g,\Gamma}"]&
\overline{{\mathcal A}}_{g,\Gamma}={\mathcal A}_{g,\Gamma}^{tor} \\
J_{g}\ar[r,phantom,"\ensuremath{\subset}"] & {\mathcal A}_{g}
\end{tikzcd}
\end{equation*}

\noindent Here $p_{g,\Gamma}\; : \;{\mathcal A}_{g,\Gamma}\longrightarrow
{\mathcal A}_{g}$ is a covering map.

\vspace{0.1in}
Finally we propose the following questions.

\vspace{0.1in}
\noindent
{\bf Question 6.1.} Let $\G$ be a neat arithmetic subgroup of $\G_g$. Does the closure
$\overline{J}_{g,\G}$ of $J_{g,\G}$ intersect the infinity boundary divisor $D_{\infty,\G}$?
If $g$ is sufficient large, it is probable that $\overline{J}_{g,\G}$ will not intersect
the boundary divisor $D_{\infty,\G}$.

\vskip 0.251cm\noindent
{\bf Question 6.2.} Let $\G$ be a neat arithmetic subgroup of $\G_g$. Does the closure
$\overline{J}_{g,\G}$ of $J_{g,\G}$ intersect the canonical divisor $K_{g,\G}$?

\vskip 0.251cm\noindent
{\bf Question 6.3.} Let $\G$ be a neat arithmetic subgroup of $\G_g$. How curved is
the closure $\overline{J}_{g,\G}$ of $J_{g,\G}$ along the boundary of $J_{g,\G}$?

\vskip 0.3cm
Quite recently using the good curvature properties of the moduli space
$(\mathcal M_g, \omega_{_{\rm WP}})$ endowed with the Weil-Petersson metric
$\omega_{_{\rm WP}}$, Liu, Sun and Yau \cite{Liu5} obtained interesting results
related to Conjecture 6.2. Let us explain their results briefly. We consider the coarse moduli space $(\mathcal M_g, \omega_{_{\rm WP}})$ endowed with the Weil-Petersson metric
$\omega_{_{\rm WP}}$ and the Siegel modular variety $(\mathcal A_g, \omega_{_{\rm H}})$ endowed with the Hodge metric $\omega_{_{\rm H}}$. Let $T_g:\mathcal M_g \lrt \mathcal A_g$ be the Torelli map (see (1.2)).
Assume that $V$ is a submanifold in $\mathcal M_g$ such that the image $T_g(V)$ is totally geodesic in
$(\mathcal A_g, \omega_{_{\rm H}})$, and also that $T_g(V)$ has finite volume. Under these two assumptions they proved that $V$ must be a ball quotient. As a corollary of this fact, it can be shown that
there is no higher rank locally symmetric subspace in $\mathcal M_g$.
A precise statement is as follows.

\begin{theorem}
Let $\Omega$ be an irreducible bounded symmetric domain and let $\Gamma\subset
{\rm Aut}(D)$ be a torsion free cocompact lattice. We set $X=\Omega/\Gamma$. Let $h$ be a canonical metric on $X$. If there exists a nonconstant holomorphic mapping
\begin{equation*}
f:(X,h)\lrt (\mathcal M_g, \omega_{_{\rm WP}}),
\end{equation*}
then $\Omega$ must be of rank $1$, i.e., $X$ must be a ball quotient.
\end{theorem}

\noindent
{\it Proof.} The proof of the above theorem can be found in \cite{Liu5}.
\hfill $\Box$

\vspace{0.2in}

\setcounter{section}{7}
\setcounter{theorem}{0}
\setcounter{equation}{0}
\renewcommand{\theequation}{7.\arabic{equation}}
\noindent{\bf 7. Final Remarks and Open Problems}
\vspace{0.1in}

In this final section we give some remarks and propose some open problems about the relations among the Schottky problem, the Andr{\'e}-Oort conjecture, Okounkov convex bodies, stable Schottky-Siegel forms, stable Schottky-Jacobi forms and the geometry of the Siegel-Jacobi space. We define the notion of stable Schottky-Jacobi forms and the concept of stable Jacobi equations for the universal hyperelliptic locus.

\vskip 0.3cm For two
positive integers $g$ and $h$, we consider the Heisenberg group
$$H_{\BR}^{(g,h)}=\big\{\,(\lambda,\mu;\ka)\,|\ \lambda,\mu\in \BR^{(h,g)},\ \kappa\in\BR^{(h,h)},\
\ka+\mu\,^t\lambda\ \text{symmetric} \big\}$$
endowed with the following multiplication law
$$\big(\lambda,\mu;\ka\big)\circ \big(\lambda',\mu';\ka'\big)=\big(\la+\la',\mu+\mu';\ka+\ka'+\la\,^t\mu'-
\mu\,^t\lambda'\big)$$
with $\big(\lambda,\mu;\ka\big),\big(\lambda',\mu';\ka'\big)\in H_{\BR}^{(g,h)}.$
We refer to \cite{YJH23, YJH24, YJH25, YJH26, YJH27, YJH16-1, YJH19, YJH22} for more details on the Heisenberg group $H_{\BR}^{(g,h)}.$
We define the {\sf Jacobi group} $G^J$ of degree $g$ and index $h$ that is the semidirect product of $Sp(2g,\BR)$ and $H_{\BR}^{(g,h)}$
$$G^J=Sp(2g,\BR)\ltimes H_{\BR}^{(g,h)}$$
endowed with the following multiplication law
$$
\big(M,(\lambda,\mu;\kappa)\big)\cdot\big(M',(\lambda',\mu';\kappa'\,)\big)
=\, \big(MM',(\tilde{\lambda}+\lambda',\tilde{\mu}+ \mu';
\kappa+\kappa'+\tilde{\lambda}\,^t\!\mu'
-\tilde{\mu}\,^t\!\lambda'\,)\big)$$
with $M,M'\in Sp(2g,\BR),
(\lambda,\mu;\kappa),\,(\lambda',\mu';\kappa') \in
H_{\BR}^{(g,h)}$ and
$(\tilde{\lambda},\tilde{\mu})=(\lambda,\mu)M'$. Then $G^J$ acts
on $\BH_g\times \BC^{(h,g)}$ transitively by
\begin{equation}
\big(M,(\lambda,\mu;\kappa)\big)\cdot
(\tau,z)=\Big( (A\tau+B)(C\tau+D)^{-1},(z+\lambda \tau+\mu)
(C\tau+D)^{-1}\Big),
\end{equation}
where $M=\begin{pmatrix} A&B\\
C&D\end{pmatrix} \in Sp(2g,\BR),\ (\lambda,\mu; \kappa)\in
H_{\BR}^{(g,h)}$ and $(\tau,z)\in \BH_g\times \BC^{(h,g)}$.

\medskip
\noindent
We note that the Jacobi group $G^J$ is {\it not} a reductive Lie group and
the homogeneous space ${\mathbb H}_g\times \BC^{(h,g)}$ is not a
symmetric space. From now on, for brevity we write
$\BH_{g,h}=\BH_g\times \BC^{(h,g)}.$ The homogeneous space $\BH_{g,h}$ is called the
{\sf Siegel-Jacobi space} of degree $g$ and index $h$.

\vskip 2mm

\vskip 0.21cm For $\tau=(\tau_{ij})\in\BH_g,$ we write $\tau=X+iY$
with $X=(x_{ij}),\ Y=(y_{ij})$ real. We put $d\tau=(d\tau_{ij})$ and $d{\overline\tau}=(d{\overline\tau}_{ij})$. We
also put
$$\PO=\,\left(\,
{ {1+\delta_{ij}}\over 2}\, { {\partial}\over {\partial \tau_{ij} }
} \,\right) \qquad\text{and}\qquad \POB=\,\left(\, {
{1+\delta_{ij}}\over 2}\, { {\partial}\over {\partial {\overline
{\tau}}_{ij} } } \,\right).$$

\vskip 0.2cm
For a coordinate $z\in\BC^{(h,g)},$ we set
\begin{eqnarray*}
z\,&=&U\,+\,iV,\quad\ \ U\,=\,(u_{kl}),\quad\ \ V\,=\,(v_{kl})\ \
\text{real},\\
dz\,&=&\,(dz_{kl}),\quad\ \ d{\overline z}=(d{\overline z}_{kl}),
\end{eqnarray*}
$$\PZ=\begin{pmatrix} {\partial}\over{\partial z_{11}} & \hdots &
 {\partial}\over{\partial z_{h1}} \\
\vdots&\ddots&\vdots\\
 {\partial}\over{\partial z_{1g}} &\hdots & {\partial}\over
{\partial z_{hg}} \end{pmatrix},\quad \PZB=\begin{pmatrix}
{\partial}\over{\partial {\overline z}_{11} }   &
\hdots&{ {\partial}\over{\partial {\overline z}_{h1} }  }\\
\vdots&\ddots&\vdots\\
{ {\partial}\over{\partial{\overline z}_{1g} }  }&\hdots &
 {\partial}\over{\partial{\overline z}_{hg} }  \end{pmatrix}.$$

\newcommand\bz{d{\overline z}}

\vskip 0.3cm
 The author proved the following theorems in \cite{YJH13}.

\begin{theorem} For any two positive real numbers
$A$ and $B$,
\begin{eqnarray}
ds_{g,h;A,B}^2&=&\,A\cdot {{\rm tr}} \Big( Y^{-1}d\tau\,Y^{-1}d{\overline
\tau}\Big) \nonumber \\
&& \ \ + \,B\cdot \bigg\{ {{\rm tr}} \Big(
Y^{-1}\,^tV\,V\,Y^{-1}d\tau\,Y^{-1} d{\overline \tau} \Big)
 +\,{{\rm tr}} \Big( Y^{-1}\,^t(dz)\,\bz\Big) \nonumber\\
&&\quad\quad -{{\rm tr}} \Big( V\,Y^{-1}d\tau\,Y^{-1}\,^t(\bz)\Big)\,
-\,{{\rm tr}} \Big( V\,Y^{-1}d{\overline \tau}\, Y^{-1}\,^t(dz)\,\Big)
\bigg\} \nonumber
\end{eqnarray}
is a Riemannian metric on $\BH_{g,h}$ which is invariant under the action {\rm (7.1)} of $G^J.$
In fact, $ds_{g,h}^2$ is a K{\"a}hler metric of $\BH_{g,h}.$
\end{theorem}
\vskip 1mm\noindent
{\it Proof.} See \cite[Theorem 1.1]{YJH13}. \hfill $\Box$

\vskip 3mm

\begin{theorem} The Laplacian $\Delta_{g,h;A,B}$ of the $G^J$-invariant metric $ds_{g,h;A,B}^2$ is given by
\begin{equation*}
\Delta_{g,h;A,B}=\,{\frac 4A}\cdot {\mathbb M}_1 + {\frac 4B}\cdot
{\mathbb M}_2,
\end{equation*}
where
\begin{eqnarray*}
{\mathbb M}_1\,&=&  {{\rm tr}} \left(\,Y\,
{}^{{}^{{}^{{}^\text{\scriptsize $t$}}}}\!\!\!
\left(Y\POB\right)\PO\,\right)\, +\,{{\rm tr}} \left(\,VY^{-1}\,^tV\,
{}^{{}^{{}^{{}^\text{\scriptsize $t$}}}}\!\!\!
\left(Y\PZB\right)\,\PZ\,\right)\\
& &\ \
+\,{{\rm tr}} \left(V\,
{}^{{}^{{}^{{}^\text{\scriptsize $t$}}}}\!\!\!
\left(Y\POB\right)\PZ\,\right)
+\,{{\rm tr}} \left(\,^tV\,
{}^{{}^{{}^{{}^\text{\scriptsize $t$}}}}\!\!\!
\left(Y\PZB\right)\PO\,\right)\nonumber
\end{eqnarray*}
and
\begin{equation*}
{\mathbb M}_2=\,{{\rm tr}}
\left(\, Y\,\PZ\,
{}^{{}^{{}^{{}^\text{\scriptsize $t$}}}}\!\!\!
\left(
\PZB\right)\,\right).
\end{equation*}
Furthermore ${\mathbb M}_1$ and ${\mathbb M}_2$ are differential operators on $\BH_{g,h}$ invariant under the action {\rm (7.1)} of $G^J.$
\end{theorem}

\noindent
{\it Proof.} See \cite[Theorem 1.2]{YJH13}. \hfill $\Box$

\vskip 3mm

\newcommand\bw{d{\overline W}}
\newcommand\bo{d{\overline \tau}}

\begin{remark}
{\rm
We refer to \cite{EXX, K-Z, YJH14, YJH16, YJH21, YY, YJH17, Z-K} for topics related 
to $ds^2_{g,h;A,B}$ and $\Box_{g,h;A,B}.$
}
\end{remark}

\begin{remark}
{\rm
Erik Balslev \cite{Bal} developed the spectral theory of $\Box_{1,1;1,1}$ on 
$\mathbb{H}_1\times\BC$ for certain arithmetic subgroups of the Jacobi modular group 
to prove that the set of all eigenvalues of $\Box_{1,1;1,1}$ satisfies the Weyl law.
}
\end{remark}

\begin{remark}
{\rm
The sectional curvature of $(\BH_1\times\BC, ds^2_{1,1;A,B})$ is $-{3\over A}$ and hence is independent of the parameter $B$. We refer to \cite{YJH17} for more detail.
}
\end{remark}

\begin{remark}
{\rm
For an application of the invariant metric $ds^2_{g,h;A,B}$ we refer to \cite{YY}.
}
\end{remark}

\begin{definition}\label{Definition 7.1}
{\rm
Let $D={\rm diag}(d_1,d_2,\cdots,d_g)$ be the $g\times g$
diagonal matrix with positive integers $d_1,\cdots,d_g$
satisfying $d_1|d_2|\cdots|d_g$,
usually called a {\sf polarization type}. $D=I_g$ is called
the {\sf principal} polarization type.
}
\end{definition}

For a fixed $\tau\in\BH_g$ and a fixed polarization type
$D={\rm diag}(d_1,\cdots,d_g)$, we let $L_\tau^D:=\BZ^g\tau+\BZ^gD$
be a lattice in $\BC^g$ and $A_\tau^D:=\BC^g/L_\tau^D$ be
a complex torus of a polarization type $D$.
Let $\{ c_0,\cdots,c_N\}$ be the set of representatives in
$\BZ^gD^{-1}$ whose components of each $c_i \,(0\leq i\leq N)$ lie in
the interval $[0,1).$ Here $N=d_1\cdots d_g-1.$

\vskip 0.251cm
We recall Lefschetz theorem (see \cite[p.\,128, Theorem 1.3]{Muf3}).

\begin{theorem}\label{Theorem 7.3}
Let $D={\rm diag}(d_1,\cdots,d_g)$ be a polarization type and $N=d_1\cdots d_g-1.$
\vskip 0.21cm\noindent
{\rm (1)} Assume $d_1\geq 2.$ Then the functions $\left\{ \theta \left[ \begin{matrix} c_0 \\ 0 \end{matrix}\right] (\tau,z),\cdots,
\theta\left[ \begin{matrix} c_N \\ 0 \end{matrix}\right] (\tau,z) \right\}$
have no zero in common, and the mapping
$\varphi^D:\BH_g\times\BC^g\lrt {\mathbb P}^N(\BC)$ defined by
\begin{equation}
\varphi^D(\tau,z):=\left( \theta \left[ \begin{matrix} c_0 \\ 0 \end{matrix}\right] (\tau,z):\cdots\cdot\cdot:
\theta\left[ \begin{matrix} c_N \\ 0 \end{matrix}\right] (\tau,z) \right),\quad (\tau,z)\in \BH_g\times\BC^g
\end{equation}
is a well-defined holomorphic mapping. For each $\tau\in\BH_g$,
the map $\varphi_\tau^D:\BC^g\lrt {\mathbb P}^N(\BC)$ defined by
\begin{equation}
\varphi^D_\tau(z):=\varphi^D(\tau,z),\qquad z\in\BC^g
\end{equation}
induces a holomorphic mapping from the complex torus $A_\tau^D$ into ${\mathbb P}^N(\BC)$.
\vskip 0.21cm\noindent
{\rm (2)} If $d_1\geq 3,$ for each $\tau\in \BH_g,$ the map
$\varphi_\tau^D:\BC^g\lrt {\mathbb P}^N(\BC)$ is an analytic embedding,
whose image is an algebraic subvariety of ${\mathbb P}^N(\BC)$.
\end{theorem}

\begin{definition}\label{Definition 7.2}
{\rm
Let $D={\rm diag}(d_1,\cdots,d_g)$ with $d_1\geq 2$ be a polarization type, and $N=d_1\cdots d_g-1.$
For each $\tau\in\BH_g,$ we define the map $\Psi^D:\BH_g\lrt {\mathbb P}^N(\BC)$ by
\begin{equation}
\Psi^D(\tau):=\varphi^D(\tau,0),\qquad \tau\in\BH_g
\end{equation}
and define the map $\Phi^D:\BH_g\times \BC^g\lrt {\mathbb P}^N(\BC)\times {\mathbb P}^N(\BC)$ by
\begin{equation}
\Phi^D(\tau,z):=\big( \varphi^D(\tau,z),\Psi^D(\tau) \big),\qquad (\tau,z)\in\BH_g\times \BC^g.
\end{equation}
}
\end{definition}

We have the following theorem proved by Baily \cite{B}.
\begin{theorem}\label{Theorem 7.4}
Assume that $d_1\geq 4$ and that $2|d_1$ or $3|d_1.$
Then the image of $\BH_g\times \BC^g$ under $\Phi^D$ is a Zariski-open subset of
an algebraic subvariety of ${\mathbb P}^N(\BC)\times {\mathbb P}^N(\BC).$
\end{theorem}
\vskip 1mm\noindent
{\it Proof.} See \cite[Section 5.1]{B} or \cite[Theorem 8.11]{P-S2}.
$\hfill \square$

\vskip 0.35cm
Let
\begin{equation*}
\G_g^J:=\G_g\ltimes H_\BZ^{(g,h)}
\end{equation*}
be the arithmetic subgroup of $G^J,$ where
\begin{equation*}
H_\BZ^{(g,h)}:=\left\{ (\lambda,\mu;\kappa)\in H_\BR^{(g,h)}\,\vert
\ \lambda,\mu,\kappa\ {\rm are\ integral}\,\right\}.
\end{equation*}
We let
\begin{equation}
{\mathcal A}_{g,h}=\G^J_g\,\ba \BH_{g,h}
\end{equation}
be the universal family of principal polarized abelian varieties of dimension $gh$. Let $\pi_{g,h}:{\mathcal A}_{g,h}\lrt
{\mathcal A}_g$
be the natural projection. We define the
{\it universal\ Jacobian\ locus}
\begin{equation}
J_{g,h}:=\pi_{g,h}^{-1}(J_g),\qquad J_g (\subset {\mathcal A}_g):={\rm the\ Jacobian\ locus}.
\end{equation}

\noindent
{\bf Problem 7.1.}
Characterize $J_{g,h}=\pi_{g,h}^{-1}(J_g).$ Describe $J_{g,h}$ in terms of Jacobi forms. We refer to
\cite{B-S, E-Z, YJH1, YJH4, YJH2, YJH6, YJH5, YJH8, YJH9, YJH10, YJH11, YJH12, YJH18, YJH21, Zi}
for more details about Jacobi forms.

\vskip 3mm\noindent
{\bf Problem 7.2.}
Compute the geodesics, the distance between two points and curvatures explicitly in the Siegel-Jacobi space $(\BH_{g,h},ds^2_{g,h;A,B}).$
See Theorem 6.1 for the Siegel space $\BH_g$.

\vskip 2mm\noindent
{\bf Problem 7.3.} Find the analogue of the Hirzebruch-Mumford Proportionality Theorem for ${\mathcal A}_{g,\G}^u$ (see (7.8) below).
\vskip 1mm
Let us give some remarks for this problem.
 Before we describe the proportionality
theorem for the Siegel modular variety, first of all we review the
compact dual of the Siegel upper half plane $\BH_g$. We note that
$\BH_g$ is biholomorphic to the generalized unit disk $\BD_g$ of
degree $g$ through the Cayley transform. We suppose that
$\Lambda=(\BZ^{2g},\langle\ ,\ \rangle)$ is a symplectic lattice with a
symplectic form $\langle\ ,\ \rangle.$ We extend scalars of the lattice
$\Lambda$ to $\BC$. Let
\begin{equation*}
{\mathfrak Y}_g:=\left\{\,L\subset \BC^{2g}\,|\ \dim_\BC L=g,\ \
\langle x,y \rangle=0\quad \textrm{for all}\ x,y\in L\,\right\}
\end{equation*}
be the complex Lagrangian Grassmannian variety parameterizing
totally isotropic subspaces of complex dimension $g$. For the
present time being, for brevity, we put $G=Sp(2g,\BR)$ and
$K=U(g).$ The complexification $G_\BC=Sp(2g,\BC)$ of $G$ acts on
${\mathfrak Y}_g$ transitively. If $H$ is the isotropy subgroup of
$G_\BC$ fixing the first summand $\BC^g$, we can identify
${\mathfrak Y}_g$ with the compact homogeneous space $G_\BC/H.$ We
let
\begin{equation*}
{\mathfrak Y}_g^+:=\big\{\,L\in {\mathfrak Y}_g\,|\ -i \langle x,{\bar
x}\rangle >0\quad \textrm{for all}\ x(\neq 0)\in L\,\big\}
\end{equation*}
be an open subset of ${\mathfrak Y}_g$. We see that $G$ acts on
${\mathfrak Y}_g^+$ transitively. It can be shown that ${\mathfrak Y}_g^+$ is biholomorphic to $G/K\cong \BH_g.$ A basis of a lattice
$L\in {\mathfrak Y}_g^+$ is given by a unique $2g\times g$ matrix
${}^t(-I_g\,\,\tau)$ with $\tau\in\BH_g$. Therefore we can identify
$L$ with $\tau$ in $\BH_g$. In this way, we embed $\BH_g$ into
${\mathfrak Y}_g$ as an open subset of ${\mathfrak Y}_g$. The
complex projective variety ${\mathfrak Y}_g$ is called the $\textit{compact dual}$ of $\BH_g.$

\newcommand\CA{\mathcal A}
\vskip 0.2cm
Let $\G$ be an arithmetic subgroup of $\G_g$. Let
$E_0$ be a $G$-equivariant holomorphic vector bundle over
$\BH_g=G/K$ of rank $r$. Then $E_0$ is defined by the
representation $\tau:K\lrt GL(r,\BC).$ That is, $E_0\cong
G\times_K \BC^r$ is a homogeneous vector bundle over $G/K$. We
naturally obtain a holomorphic vector bundle $E$ over
$\CA_{g,\G}:=\G\ba G/K.$ $E$ is often called an
$\textit{automorphic}$ or $ \textit{arithmetic}$ vector bundle
over $\CA_{g,\G}$. Since $K$ is compact, $E_0$ carries a
$G$-equivariant Hermitian metric $h_0$ which induces a Hermitian
metric $h$ on $E$. According to Main Theorem in \cite{Muf2}, $E$
admits a $ \textit{unique}$ extension ${\tilde E}$ to a smooth
toroidal compactification ${\tilde \CA}_{g,\G}$ of $\CA_{g,\G}$
such that $h$ is a singular Hermitian metric $ \textit{good}$ on
${\tilde \CA}_{g,\G}$. For the precise definition of a
$\textit{good metric}$ on $\CA_{g,\G}$ we refer to \cite[p.\,242]{Muf2}.
According to Hirzebruch-Mumford's Proportionality
Theorem\,(cf.~\cite[p.\,262]{Muf2}), there is a natural metric on
$G/K=\BH_g$ such that the Chern numbers satisfy the following
relation
\begin{equation*}
c^{\al}\big({\tilde E}\big)=(-1)^{{\frac 12}g(g+1)}\,
\textmd{vol}\left( \G\ba \BH_g\right)\,c^{\al}\big( {\check E}_0\big)
\end{equation*}
for all $\al=(\al_1,\cdots,\al_r)$ with nonegative integers
$\al_i\,(1\leq i\leq r)$ and $\sum_{i=1}^r\al_i={\frac 12}g(g+1),$
where ${\check E}_0$ is the $G_{\BC}$-equivariant holomorphic
vector bundle on the compact dual ${\mathfrak Y}_g$ of $\BH_g$
defined by a certain representation of the stabilizer $
\textrm{Stab}_{G_\BC}(e)$ of a point $e$ in ${\mathfrak Y}_g$.
Here $\textmd{vol}\left( \G\ba \BH_g\right)$ is the volume of
$\G\ba\BH_g$ that can be computed\,(cf.~\cite{Si1}).

\vskip 3mm\noindent
{\bf Problem 7.4.} Compute the cohomology
$H^\bullet ({\mathcal A}_{g,h},*)$ of ${\mathcal A}_{g,h}.$
Investigate the intersection cohomology of ${\mathcal A}_{g,h}.$

\vskip 2mm\noindent
{\bf Problem 7.5.} Generalize the trace formula on the Siegel modular variety obtained by Sophie Morel to the universal abelian variety. For her result on the trace formula on the Siegel modular variety, we refer to her paper \cite[{\it Cohomologie d'intersection des vari{\'e}t{\'e}s modulaires de Siegel, suite}]{Mor}.

\vskip 3mm\noindent
{\bf Problem 7.6.}
Construct all the geodesics contained in $J_{g,h}.$

\vskip 3mm\noindent
{\bf Problem 7.7.}
Develop the theory of variations of abelian varieties along the geodesic joining two points in $J_{g,h}$.

\vskip 3mm\noindent
{\bf Problem 7.8.}
Discuss the Andr{\'e}-Oort conjecture for ${\mathcal A}_{g,h}$. Gao proved the Ax-Lindemann-Weierstras theorem for ${\mathcal A}_{g,h}$, and using this theorem proved
the Andr{\'e}-Oort conjecture for ${\mathcal A}_{g,h}$
under the assumption of the Generalized Riemann Hypothesis for CM fields in his paper \cite{G}.

\vskip 3mm
Let $\G$ be a neat arithmetic subgroup of $\G_g$. We put $\G^J:=\G\ltimes H_\BZ^{(g,h)}$. We let
\begin{equation}
{\mathcal A}_{g,h,\G}:=\G^J\,\backslash \BH_{g,h}.
\end{equation}
Let ${\mathcal A}_{g,h,\G}^{tor}$ be a toroidal compactification of ${\mathcal A}_{g,h,\G}$. Let
$K_{g,h,\G}$ be the canonical line bundle over
${\mathcal A}_{g,h,\G}^{tor}$ and let
\begin{equation}
D_{\infty,g,h,\G}:={\mathcal A}_{g,h,\G}^{tor}\backslash
{\mathcal A}_{g,h,\G}
\end{equation}
be the infinity boundary divisor on ${\mathcal A}_{g,h,\G}^{tor}$. Let $\pi_{g,h,\G}:{\mathcal A}_{g,h,\G}\lrt {\mathcal A}_{g,\G}$
be a projection and let $p_{g,\G}:{\mathcal A}_{g,\G}\lrt
{\mathcal A}_{g}$ be a covering map. We define
\begin{equation}
J_{g,h,\G}:=\left( p_{g,\G}\circ \pi_{g,h,\G}\right)^{-1}(J_g).
\end{equation}

\noindent
{\bf Problem 7.9.} Assume that ${\mathcal A}_{4,h,\G}^{tor}$ is a toroidal compactification of ${\mathcal A}_{4,h,\G}$ which is projective.
Compute the Okounkov bodies $\Delta_{Y_\bullet}(J_{4,h,\G})$, $\Delta_{Y_\bullet}^{\rm val}(J_{4,h,\G})$ and
$\Delta_{Y_\bullet}^{\rm lim}(J_{4,h,\G})$ explicitly. Describe the relations among $J_4,\ J_{4,h,\G}$ and these Okounkov bodies explicitly.
Describe the relations between these Okounkov bodies and the $GL(4,\BZ)$-admissible family of polyhedral decompositions defining the toroidal compactification $\overline{\mathcal A}_{4,\G}$.

\vskip 3mm\noindent
{\bf Problem 7.10.} Assume that a toroidal compactification ${\mathcal A}_{g,h,\G}^{tor}$ is a projective
variety. Let $K_{g,h,\G}$ be the canonical line bundle over ${\mathcal A}_{g,h,\G}^{tor}$ and
$D_{\infty,g,h,\G}$ be the infinity boundary divisor on
${\mathcal A}_{g,h,\G}^{tor}$.
Compute the Okoukov convex bodies $\Delta_{Y_\bullet}(K_{g,h,\G}),\ 
\Delta_{Y_\bullet}^{\rm val}(K_{g,h,\G}),
\Delta_{Y_\bullet}^{\rm lim}(K_{g,h,\G}), \Delta_{Y_\bullet}(D_{\infty,g,h,\G}),
\Delta_{Y_\bullet}^{\rm val}(D_{\infty,\G}^u),
\ \Delta_{Y_\bullet}^{\rm lim}(D_{\infty,g,h,\G}),\\ 
\Delta_{Y_\bullet}(K_{g,h,\G}+D_{\infty,g,h,\G}),
\ \Delta_{Y_\bullet}^{\rm val}(K_{g,h,\G}+D_{\infty,g,h,\G})$ and
$\Delta_{Y_\bullet}^{\rm lim}(K_{g,h,\G}+D_{\infty,g,h,\G})$ explicitly. Describe the relations between these Okounkov bodies and the $GL(g,\BZ)$-admissible family of polyhedral decompositions defining the toroidal compactification
${\mathcal A}_{g,h,\G}^{tor}$.

\vskip 3mm\noindent
{\bf Problem 7.11.} Assume that a toroidal compactification ${\mathcal A}_{g,h,\G}^{tor}$ of ${\mathcal A}_{g,h,\G}$ is a projective variety.
Let $D_{J,\G}$ be a divisor on ${\mathcal A}_{g,h,\G}^{tor}$ containing $J_{g,h,\G}$. Describe the Okounkov bodies
$\Delta_{Y_\bullet}(D_{J,\G}),\ \Delta_{Y_\bullet}^{\rm val}(D_{J,\G})$ and $\Delta_{Y_\bullet}^{\rm lim}(D_{J,\G})$.
Study the relations among $J_{g,\G},\ J_{g,h,\G},\ D_{J,\G}$ and these Okounkov bodies.

\vskip 0.251cm
We have the following diagram:
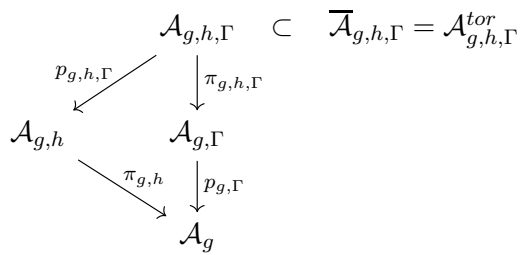
\begin{equation*}
\begin{tikzcd}[ampersand replacement = \&]
\& {\mathcal A}_{g,h,\Gamma}  \ar[r,phantom,"\ensuremath{\subset}"] \ar[d,"\pi_{g,h,\Gamma}"]\ar[dl,outer sep=-2pt, "p_{g,h,\Gamma} "']
\&  \overline{{\mathcal A}}_{g,h,\Gamma}
={\mathcal A}^{tor}_{g,h,\Gamma}  \\
{\mathcal A}_{g,h}  \ar[dr,outer sep=-2pt, "\pi_{g,h}"]
\&{\mathcal A}_{g,\Gamma}\ar[d,"p_{g,\Gamma}"]  \\
\&  {\mathcal A}_{g}
\end{tikzcd}
\end{equation*}

\noindent Here $p_{g,h,\Gamma}\; : \;
{\mathcal A}_{g,h,\Gamma}\longrightarrow {\mathcal A}_{g,h}$ is a covering map.

\vskip 0.35cm
We propose the following questions.
\vskip 0.251cm\noindent
{\bf Question 7.1.} Let $\G$ be a neat arithmetic subgroup of $\G_g$. Does the closure
$\overline{J}_{g,h,\G}$ of $J_{g,h,\G}$ intersect the infinity boundary divisor $D_{\infty,g,h,\G}$?
If $g$ is sufficient large, it is probable that $\overline{J}_{g,h,\G}$ will not intersect
the boundary divisor $D_{\infty,g,h,\G}^u$.

\vskip 0.251cm\noindent
{\bf Question 7.2.} Let $\G$ be a neat arithmetic subgroup of $\G_g$. Does the closure
$\overline{J}_{g,h,\G}$ of $J_{g,h,\G}$ intersect the canonical divisor $K_{g,h,\G}$?

\vskip 0.251cm\noindent
{\bf Question 7.3.} Let $\G$ be a neat arithmetic subgroup of $\G_g$. How curved is
the closure $\overline{J}_{g,h,\G}$ of $J_{g,h,\G}$ along the boundary of $J_{g,h,\G}$?

\vskip 0.35cm Now we make some conjectures.

\begin{conjecture}
For a sufficiently large integer $g$, the locus $J_{g,h}$ contains only finitely many special points. This is an analogue $($or generalization$)$ of Coleman's conjecture.
\end{conjecture}

\begin{conjecture}
For a sufficiently large integer $g$, the locus $J_{g,h}$ cannot contain a non-trivial totally geodesic subvariety inside
${\mathcal A}_{g,h}$ for the Riemannian metric $ds^2_{g,h;A,B}$.
\end{conjecture}

\begin{conjecture}
For a sufficiently large integer $g$, there does not exist a geodesic that is contained in $J_{g,h}$
for the Riemannian metric $ds^2_{g,h;A,B}$.
\end{conjecture}

Finally we discuss the connection between the universal Jacobian locus $J_{g,h}$ and
stable Jacobi forms. We refer to Appendix E in this article for more details on stable Jacobi forms. First we review the concept of {\it stable\ modular\ forms} introduced in
\cite{Fr2}. The Siegel $\Phi$-operator
\begin{equation}
\Phi_{g,k}:[\G_{g+1},k]\lrt [\G_g,k], \qquad k\in\BZ_+
\end{equation}
defined by
\begin{equation*}
(\Phi_{g,k}f)(\tau):=\lim_{t\rightarrow \infty}\,f\left(
\begin{pmatrix} \tau & 0\\ 0 & it\end{pmatrix}\right), \quad f\in [\G_{g+1},k],
\,\tau\in \BH_g,
\end{equation*}
where $[\G_g,k]$ denotes the vector space of all Siegel modular forms on $\BH_g$ of weight $k$. Using the theory of Poincar{\'e} series, H. Maass \cite{M1} proved that
if $k$ is even and $k>2g$, then $\Phi_{g,k}$ is a surjective linear map. In 1977,
using the theory of singular modular forms, E. Freitag \cite{Fr2} proved the following facts (a) and (b)~:
\vskip 1mm
(a) for a fixed even integer $k$, $\Phi_{g,k}$ is an isomorphism if $g>2k$~;
\vskip 1mm
(b) $[\Gamma_g,k]=0$\quad if $g>2k,\ k\not\equiv 0~({\rm mod}~4).$\emph{}

\vskip 3mm
The fact (a) means that the vector spaces $[\G_g,k]$ stabilize to the infinity vector space $[\G_{\infty}, k]$ as $g$ increases. In this sense, he introduced the notion of the {\sf stability} of Siegel modular forms.

\begin{definition}
{\rm
A Siegel modular form $f\in [\Gamma_g,k]$ is said to be {\sf stable} if there exists a nonegative integer $m\in \BZ_+$ satisfying the following conditions {\rm (SM1)} and {\rm (SM2)}~:
\vspace{0.1in}

{\rm (SM1)} $g+m> 2k$;
\vskip 1mm
{\rm (SM2)} $f=\Phi_{g+1,k}\circ \Phi_{g+2,k}\circ\cdots\circ\Phi_{g+m,k}(F)\quad$ for some
$F\in [\G_{g+m},k]$.
}
\end{definition}

Scalar-valued Siegel modular forms on $\mathcal A_g$ vanishing on the Jacobian locus, equivalently, forms on the Satake compactification
$\mathcal A_g^{\rm Sat}$ that vanish on the closure $J_g^{\rm Sat}$ of $J_g$ in $\mathcal A_g^{\rm Sat}$ are called {\sf Schottky-Siegel forms}. The normalization $\nu:\mathcal A_g^{\rm Sat}\lrt \partial\mathcal A_{g+1}^{\rm Sat}$ gives a restriction map which coincides with the Siegel operator $\Phi_{g,k}~(k\in\BZ_+).$

\vskip 2mm
We let
$$A(\G_g):=\bigoplus_{k\geq 0} [\G_g,k]$$
be the graded ring of Siegel modular forms on $\BH_g.$ It is known that $A(\G_g)$ is a finitely generated $\BC$-algebra and the field of modular functions $K(\Gamma_g)$ is an algebraic function field of transcendence degree ${\frac 12}g(g+1)$.

The ring
\begin{equation*}
{\mathbb A}=\bigoplus_{k\geq 0} [\G_\infty,k]
\end{equation*}
is an inverse limit in the category
\begin{equation}
{\mathbb A}=\lim_{\begin{subarray}{c} \longleftarrow\\ ^g \end{subarray}} A(\G_g).
\end{equation}
Freitag \cite{Fr2} proved that ${\mathbb A}$ is the polynomial ring over $\BC$ on the set of theta series $\theta_S$, where $S$ runs over the set of equivalence classes of indecomposable positive definite unimodular even integral matrices. In general, $A(\G_g)$ is not a polynomial ring (cf.~\cite[p.\,204]{Fr2}).

\vskip 2mm
We define the {\sf stable} Satake compactification $\mathcal A_\infty^{\rm Sat}$ by
\begin{equation}
\mathcal A_\infty^{\rm Sat}:=\bigcup_g \mathcal A_g^{\rm Sat}=
\lim_{\begin{subarray}{c} \longleftarrow\\ ^g \end{subarray}} \mathcal A_g^{\rm Sat}
\end{equation}
and the {\sf stable Jacobian locus} $J_\infty^{\rm Sat}$ by
\begin{equation}
J_\infty^{\rm Sat}:=\bigcup_g J_g^{\rm Sat}=
\lim_{\begin{subarray}{c} \longleftarrow\\ ^g \end{subarray}} J_g^{\rm Sat}.
\end{equation}

\vskip 2mm
G. Codogni and N.~I. Shepherd-Barron \cite{Cod-Sh} proved the following theorem.

\begin{theorem}
There are no stable Schottky-Siegel forms. That is, the homomorphism from
\begin{equation}
{\mathbb A}=\lim_{\begin{subarray}{c} \longleftarrow\\ ^g \end{subarray}}
A(\G_g)\lrt \bigoplus_k H^0
(J_\infty^{\rm Sat},
\omega_{_{J}}^{\otimes k})
\end{equation}
induced by the inclusion $J_\infty^{\rm Sat} \hookrightarrow
\mathcal A_\infty^{\rm Sat}$ is injective, where $\omega_{_{J}}$ is the restriction of the canonical line bundle $\omega$ on $\mathcal A_\infty^{\rm Sat}$ to
$J_\infty^{\rm Sat}$.
\end{theorem}

\noindent
{\it Proof.} See Theorem 1.3 and Corollary 1.4 in \cite{Cod-Sh}. \hfill $\Box$

\vskip 3mm
We refer to Appendix D in this paper for the definition of Jacobi forms.
\vskip 2mm
\vskip 2mm Now we consider the special case $\rho=\det^k$ with $k\in \BZ_+$. We define the {\sf Siegel-Jacobi operator}
\begin{equation*}
\Psi_{g,\M}: J_{k,\CM}(\G_g)\lrt J_{k,\CM}(\G_{g-1})
\end{equation*}
by
\begin{equation}
(\Psi_{g,\M}F)(\tau,z):= \lim_{t\lrt\infty}
F\left( \begin{pmatrix} \tau & 0 \\ 0 & it \end{pmatrix},(z,0)\right),\\
\end{equation}
where $F \in J_{k,\CM}(\G_g),\ \tau\in \BH_{g-1}$ and
$z\in \BC^{(h,g-1)}.$ We observe that the above limit exists
and $\Psi_{g,\M}$ is a well-defined linear map\,(cf.\,\cite{Zi}).

\vspace{0.1in}
The author\,\cite{YJH4} proved the following theorems.
\begin{theorem}
Let $2\M$ be a positive even unimodular symmetric integral matrix of degree $h$ and let $k$ be an even nonnegative integer. If $g+h>2k$, then the Siegel-Jacobi operator
$\Psi_{g,\M}$ is injective.
\end{theorem}

\noindent
{\it Proof.} See \cite[Theorem 3.5]{YJH4}. \hfill $\Box$

\begin{theorem}
Let $2\M$ be as above in {\rm Theorem 2.1} and let $k$ be an even nonnegative integer. If $g+h>2k+1$, then the Siegel-Jacobi operator $\Psi_{g,\M}$ is an isomorphism.
\end{theorem}
%
\noindent
{\it Proof.} See \cite[Theorem 3.6]{YJH4}. \hfill $\Box$

\begin{theorem}
Let $2\M$ be as above in {\rm Theorem 2.1} and let $k$ be an even nonnegative integer. Assume that $2k>4g+h$ and
$k\equiv 0\,({\rm mod}\,2).$ Then the Siegel-Jacobi operator $\Psi_{g,\M}$ is surjective.
\end{theorem}

\noindent
{\it Proof.} See \cite[Theorem 3.7]{YJH4}. \hfill $\Box$

\begin{remark}
{\rm
 The author \cite[Theorem 4.2]{YJH4} proved that the action of the Hecke operatos on Jacobi forms is compatible with that of the Siegel-Jacobi operator.
}
\end{remark}

\begin{definition}
{\rm
A collection $(F_g)_{g\geq 0}$ is called a {\sf stable Jacobi form} of weight $k$ and index $\M$ if it satisfies the following conditions {\rm (SJ1)} and {\rm (SJ2)}:
\vskip 2mm
{\rm (SJ1)} \ \ $F_g\in J_{k,\CM}(\G_g)$ \ \ for all $g\geq 0.$
\vskip 2mm
{\rm (SJ2)} \ \ $\Psi_{g,\M}F_g=F_{g-1}$ \ \ for all $g\geq 1.$
}
\end{definition}

\begin{remark}
{\rm
The concept of a stable Jacobi forms was introduced by the author\,\cite{YJH3, YJH9}.
}
\end{remark}

\noindent {\bf Example.}
{\rm
Let $S$ be a positive even unimodular symmetric integral matrix of degree $2k$ and let $c\in \BZ^{(2k,h)}$ be an integral matrix. We define the theta series $\vartheta_{S,c}^{(g)}$ by
\begin{equation*}
\vartheta_{S,c}^{(g)}(\tau,z):=\sum_{\la\in\BZ^{(2k,g)}}
e^{\pi i\left\{ {{\rm tr}} (S\la\tau\,{}^t\la)+ 2\,{{\rm tr}} ({}^tc\,S\la\,{}^tz)\right\} }, \quad (\tau,z)\in \BH_{g,h}.
\end{equation*}
It is easily seen that $\vartheta_{S,c}^{(g)}\in J_{k,\CM}(\G_g)$
with $\M:={\frac 12}~\!{}^t\!cSc$ for all $g\geq 0$ and $\Psi_{g,\M}\vartheta_{S,c}^{(g)}=\vartheta_{S,c}^{(g-1)}$ for all
$g\geq 1.$ Thus the collection
\begin{equation*}
\Theta_{S,c}:=\left( \vartheta_{S,c}^{(g)} \right)_{g\geq 0}
\end{equation*}
is a stable Jacobi form of weight $k$ and index $\M$.
}

\begin{definition}
{\rm
Let $\M$ be a half-integral semi-positive symmetric matrix of degree $h$ and $k\in\BZ_+.$ A Jacobi form $F\in J_{k,\CM}(\G_g)$ is called a {\sf Schottky-Jacobi form} of weight $k$ and index $\M$ for the universal Jacobian locus if it vanishes along $J_{g,h}$.
}
\end{definition}

\begin{definition}
{\rm
Let $\M$ be a half-integral semi-positive symmetric matrix of degree $h$ and $k\in\BZ_+.$
A collection $(F_g)_{g\geq 0}$ is called a {\sf stable Schottky-Jacobi form} of weight $k$ and index $\M$ if it satisfies the following conditions {\rm (1)} and {\rm (2)}:
\vskip 2mm
{\rm (1)} \ $F_g\in J_{k,\M}(\G_g)$ is a Schottky-Jacobi form of weight $k$ and index $\M$ for all\\
 \indent \ \ \ \ \ \ $g\geq 0.$
\vskip 2mm
{\rm (2)}\ $\Psi_{g,\M}F_g=F_{g-1}$ for all $g\geq 1.$
}
\end{definition}

We expect to prove the following claim~:
\vskip 2mm\noindent
{\bf Claim:} {\it There are no stable Schottky-Jacobi forms for the universal Jacoban locus.}

\vskip 3mm The author \cite{YJH28} proved the following.

\begin{theorem}
Let $2\M$ be a positive even unimodular symmetric integral matrix of degree $h$. Then there do not exist stable Schottky-Jacobi forms of index $\M$ for the universal Jacobian locus.
\end{theorem}

\noindent
{\it Proof.} See \cite[Theorem 4.1]{YJH28}. \hfill $\Box$

\vskip 3mm
Let $(\Lambda, Q)$ be an even unimodular positive definite quadratic form of rank $m$. That is, $\Lambda$ is a finitely generated free group of rank $m$ and $Q$ is an integer-valued bilinear form on $\Lambda$ such that $Q$ is even and unimodular.
For a positive integer $g$, the theta series $\theta_{Q,g}$ associated to $(\Lambda, Q)$ is defined to be
\begin{equation*}
\theta_{Q,g}(\tau):=\sum_{x_1,\cdots,x_g\in \Lambda}\exp \left( \pi i \sum_{p,q}^g Q(x_p,x_q)\tau_{pq}\right),\qquad \tau=(\tau_{pq})\in \BH_g.
\end{equation*}
It is well known that $\theta_{Q,g}(\tau)$ is a Siegel modular form on $\BH_g$ of weight ${\frac m2}.$ We easily see that
\begin{equation*}
\Phi_{g,{\frac m2}}(\theta_{Q,g+1})=\theta_{Q,g}\qquad {\rm for\ all}\ g\in \BZ_+.
\end{equation*}
Therefore the collection of all theta series associated to $(\Lambda, Q)$
\begin{equation*}
\Theta_Q:=\left( \theta_{Q,g} \right)_{g\geq 0}
\end{equation*}
is a stable modular form.

\begin{definition}
{\rm
A {\sf stable equation} for the hyperelliptic locus is a stable modular form
$(f_g)_{g\geq 0}$ such that $f_g$ vanishes along the hyperelliptic locus ${\rm Hyp}_g$ for every $g$.
}
\end{definition}

Recently G. Codogni \cite{Cod} proved the following.
\begin{theorem}
The ideal of stable equations of the hyperelliptic locus is generated by differences of theta series
$$\theta_P-\theta_Q,$$
where $P$ and $Q$ are even unimodular positive definite quadratic forms of the same rank.
\end{theorem}
\vskip 1mm\noindent
{\it Proof.} See Theorem 1.2 or Theorem 4.2 in \cite{Cod}. \hfill $\Box$

\vskip 3mm
In a similar way we may define the concept of stable Jacobi equation.

\begin{definition}
{\rm
A {\sf stable Jacobi equation} of index $\mathcal M$ for the universal hyperelliptic locus is a stable Jacobi form $(F_{g,\M})_{g\geq 0}$ of index $\mathcal M$ such that $F_{g,\M}$ vanishes along the universal hyperelliptic locus ${\rm Hyp}_{g,h}:=\pi_{g,h}^{-1}({\rm Hyp}_g)$ for every $g$.
}
\end{definition}

The author \cite{YJH28} proved the following.
\begin{theorem}
Let $2\M$ be a positive even unimodular symmetric integral matrix of degree $h$. Then there exist non-trivial stable Schottky-Jacobi forms of $\M$ for the universal hyperelliptic locus.
\end{theorem}

\noindent
{\it Proof.} See \cite[Theorem 4.2]{YJH28}. \hfill $\Box$

\vskip 3mm
\noindent
{\bf Problem 7.12.} Find the ideal of stable Jacobi equations of the universal hyperelliptic locus.

\begin{remark}
{\rm
We consider a half-integral semi-positive symmetric integral
matrix $\M$ such that $2\M$ is {\it not} even or which is {\it not} unimodular.
The natural questions arise:
\vskip 2mm\noindent
{\bf Question\ 7.1.} Are there non-trivial stable Schottky-Jacobi forms of index $\M$ for the universal Jacobian locus?
\vskip 2mm\noindent
{\bf Question\ 7.2.} Are there non-trivial stable Schottky-Jacobi forms of index $\M$ for the universal hyperelliptic locus?
}
\end{remark}

\vspace{0.1in}


\def\BP{\Bbb P}
\def\vp{\varphi}
\def\J{J\in {\Bbb Z}^{(h,g)}_{\geq 0}}
\def\N{N\in {\Bbb Z}^{(h,g)}}
\def\wt{\widetilde}
\def\Box{$\square$}
\def\M{\mathcal M}
\def\ta{\vartheta^{(\mathcal{M})}\left[\begin{pmatrix} A\\ 0\end{pmatrix}\right]
(\Omega\vert W)}
\def\Jt{\left({{\partial}\over {\partial W}}\right)^J\vartheta^{(\mathcal{M})}
\left[\begin{pmatrix} A\\ 0\end{pmatrix}\right](\Omega\vert W)}
\def\Dt{\Delta^J\vartheta^{(\mathcal{M})}\left[\begin{pmatrix} A\\ 0\end{pmatrix}\right]
(\Omega\vert W)}
\def\emw{\text {exp}\left\{-\pi i\sigma({\mathcal{M}}(\xi\Omega\,^t\xi
+2W\,^t\!\xi))\right\} }
\def\mjz{{\tilde {\vartheta}}^{(\mathcal{M})}_J\left[\begin{pmatrix} A\\ 0\end{pmatrix}
\right](\Omega\vert Z,W)}
\def\CJ{{\Bbb C}\left[\cdots,\left({{\partial}\over
{\partial W}}\right)^J
\vartheta^{(\mathcal{M})}\left[\begin{pmatrix} A_{\alpha,{\mathcal{M}}}\\ 0\end{pmatrix}
\right](\Omega\vert W),\cdots\right]}
\def\GJ{G\left(\cdots,\Delta^J\vartheta^{(\mathcal{M})}\left[\begin{pmatrix}
A_{\alpha,{\mathcal{M}}}\\ 0\end{pmatrix}\right](\Omega\vert W),\cdots\right)}
\def\GW{G\left(\cdots,\left({{\partial}\over {\partial W}}\right)^J
\vartheta^{(\mathcal{M})}\left[\begin{pmatrix} A_{\alpha,{\mathcal{M}}}\\ 0\end{pmatrix}\right]
(\Omega\vert W),\cdots\right)}
\def\lb{\lbrace}
\def\lk{\lbrack}
\def\rb{\rbrace}
\def\rk{\rbrack}
\def\s{\sigma}
\def\w{\wedge}
\def\lrt{\longrightarrow}
\def\lmt{\longmapsto}
\def\O{\Omega}
\def\k{\kappa}
\def\ba{\backslash}
\def\ph{\phi}
\def\H{\mathbb H}
\def\l{\lambda}
\def\A{\mathcal A}

\setcounter{theorem}{0}
\setcounter{equation}{0}
\noindent{\bf Appendix A. Subvarieties of the Siegel Modular Variety}
\vspace{0.1in}

In this appendix A, we give a brief remark on subvarieties of the Siegel modular variety and present several problems. This appendix was written on the base of the review
\cite{San} of G.~K. Sankaran for the paper \cite{YJH15}. In fact, Sankaran made a critical review on Section 10. Subvarieties of the Siegel modular variety of the
author's paper \cite{YJH15} and corrected some wrong statements and information given
by the author. In this sense the author would like to give his deep thanks to the reviewer, Sankanran.

\vspace{0.1in}
\indent Here we assume that the ground field is the complex number field $\BC.$

\vspace{0.1in}
\noindent{\bf Definition A.1.} A nonsingular variety $X$ is said to be
{\it rational} if $X$ is birational to a projective space $\BP^n(\BC)$ for
some integer $n$. A nonsingular variety $X$ is said
to be {\it stably\ rational} if $X\times \BP^k(\BC)$ is birational to $\BP^N(\BC)$ for certain nonnegative
integers $k$ and $N$. A nonsingular variety $X$ is called {\it unirational} if there
exists a dominant rational map $\varphi:\BP^n(\BC)\lrt X$ for a certain
positive integer $n$, equivalently if the function field $\BC(X)$ of $X$
can be embedded in a purely transcendental extension $\BC(z_1,\cdots,z_n)$ of $\BC.$
\vspace{0.1in}\\
\noindent{\bf Remarks A.2.} (1) It is easy to see that the rationality implies
the stably rationality and that the stably rationality implies the
unirationality.
\vspace{0.05in}\\
\noindent (2) If $X$ is a Riemann surface or a complex surface, then the notions of
rationality, stably rationality and unirationality are equivalent one
another.
\vspace{0.05in}\\
\noindent (3) H. Clemens and P. Griffiths \cite{Cl-G} showed that most of cubic threefolds in $\BP^4(\BC)$ are unirational but {\it not} rational.
\vspace{0.1in}\\
\indent
The following natural questions arise :
\vspace{0.1in}\\
{\bf Question 1.} Is a stably rational variety {\it rational} ?
\vspace{0.1in}\\
{\bf Question 2.} Is a general hypersurface $X\subset \BP^{n+1}(\BC)$ of
degree $d\leq n+1$ {\it unirational} ?

\vskip 3mm
Question 1 is a famous one raised by O. Zariski (cf.~B. Serge, {\sf Algebra and Number Theory (French)}, CNRS, Paris (1950), 135--138; MR0041480). In \cite{BCSS}, A. Beauville, J.-L. Colliot-Th{\'e}l{\`e}ne, J.-J. Sansuc and P. Swinnerton-Dyer
gave counterexamples, e.g., the Ch${\hat \rm a}$telot surfaces $V_{d,P}\subset {\bf A}_\BC^3$ defined by $y^2-dz^2=P(x),$ where $P\in\BC [x]$ is an irreducible polynomial of degree 3, and $d$ is the discriminant of $P$ such that $d$ is not a square and hence answered negatively to Question 1.
\vspace{0.1in}\\
\noindent{\bf Definition A.3.} Let $X$ be a nonsingular variety of dimension $n$
and let $K_X$ be the canonical divisor of $X$. For each positive integer
$m\in \BZ^+$, we define the $m$-{\it genus} $P_m(X)$ of $X$ by
$$P_m(X):=\text{dim}_{\BC}\,H^0(X,{\mathcal{O}}(mK_X)).$$
The number $p_g(X):=P_1(X)$ is called the {\it geometric\ genus} of $X$.
We let
$$N(X):=\left\{\,m\in \BZ^+\,\vert\, P_m(X)\geq 1 \right\}.$$
For the present, we assume that $N(X)$ is nonempty. For each $m\in N(X),$
we let $\left\{ \phi_0,\cdots,\phi_{N_m}\right\}$ be a basis of the vector
space $H^0(X,{\mathcal{O}}(mK_X)).$ Then we have the mapping $\Phi_{mK_X}\,:
X\lrt \BP^{N_m}(\BC)$ by
$$\Phi_{mK_X}(z):=(\phi_0(z):\cdots:\phi_{N^m}(z)),\ \ \ z\in X.$$
We define the {\it Kodaira\ dimension} $\kappa(X)$ of $X$ by
$$\kappa(X):=\text{max}\,\left\{\,\text{dim}_{\BC}\,\Phi_{mK_X}(X)\,\vert\,\,
m\in N(X) \right\}.$$
If $N(X)$ is empty, we put $\kappa(X):=-\infty.$ Obviously $\kappa(X)\leq
\text{dim}_{\BC}\,X.$ A nonsingular variety $X$
is said to be {\it of\ general\
type} if $\kappa(X)=\text{dim}_{\BC}X.$
A singular variety $Y$ in general is said to be rational, stably rational,
unirational or of general type if any nonsingular model $X$ of $Y$ is
rational, stably rational, unirational or of general type respectively.
We define
$$P_m(Y):=P_m(X)\ \ \text{and}\ \ \kappa(Y):=\kappa(X).$$
A variety $Y$ of dimension $n$ is said to be {\it of\ logarithmic\
general\ type} if there exists a smooth compactification ${\tilde Y}$ of
$Y$ such that $D:={\tilde Y}-Y$ is a divisor with normal crossings only and
the transcendence degree of the logarithmic canonical ring
$$\oplus_{m=0}^{\infty}\,H^0({\tilde Y},\,m(K_{\tilde Y}+[D]))$$
is $n+1$, i.e., the {\it logarithmic\ Kodaira\ dimension} of $Y$ is $n$.
We observe that the notion of being of logarithmic general type is weaker
than that of being of general type.
\vspace{0.1in}\\
\indent Let $\A_g:=\G_g\ba \H_g$ be the Siegel modular variety of degree $g$, that is,
the moduli space of principally polarized abelian varieties of dimension
$g$. So far it has been proved that $\A_g$ is of general type for $g\geq 7.$
At first Freitag \cite{Fr1} proved this fact when $g$ is a multiple of $24$.
Tai \cite{T} proved this for $g\geq 9$ and Mumford \cite{Muf2-1} proved this fact for
$g\geq 7.$ On the other hand, $\A_g$ is known to be unirational for
$g\leq 5\,:$ Donagi \cite{DO0} for $g=5,$ Clemens \cite{Cl} for $g=4$ and classical for
$g\leq 3.$ For $g=3,$ using the moduli theory of curves, Riemann \cite{Rie},
Weber \cite{We} and Frobenius \cite{Frob} showed that $\A_3(2):=\G_3(2)\ba \H_3$ is a
rational variety and moreover gave $6$ generators of the modular function
field $K(\G_3(2))$ written explicitly in terms of derivatives of odd theta
functions at the origin. So $\A_3$ is a unirational variety with a Galois
covering of a rational variety of degree $[\G_3:\G_3(2)]=1,451,520.$ Here
$\G_3(2)$ denotes the principal congruence subgroup of $\G_3$ of level $2.$
Furthermore it was shown that $\A_3$ is stably rational~(cf. \cite{Bo-K,KS}).
For a positive integer $k$, we let $\G_g(k)$ be the principal
congruence subgroup of $\G_g$ of level $k$. Let $\A_g(k)$ be the moduli
space of abelian varieties of dimension $g$ with $k$-level structure. It is
classically known that $\A_g(k)$ is of logarithmic general type for
$k\geq 3$ (cf. \cite{Muf2-1}). Wang \cite{W1, W2} gave a different proof for the fact that $\A_2(k)$ is of general type for $k\geq 4.$ On the other hand, the relation between the Burkhardt quartic and abelian surfaces with 3-level structure was established by H. Burkhardt \cite{Bu} in 1890. We refer to \cite[\S\, IV.2,~pp.\,132--135]{H-S} for more detail on the Burkhardt quartic. In 1936, J.~A. Todd \cite{To} proved that the Burkhardt quartic is rational.
van der Geer \cite{vdG2} gave a modern proof for the rationality of $\A_2(3)$.
The remaining unsolved problems are summarized as follows:

\vspace{0.1in}
\noindent
{\bf Problem 1.} Are $\A_4,\ \A_5$ stably rational or rational?

\vspace{0.1in}
\noindent
{\bf Problem 2.} Discuss the (uni)rationality of $\A_6.$

\vspace{0.1in}
\indent We already mentioned that $\A_g$ is of general type if $g\geq 7.$ It is
natural to ask if the subvarieties of $\A_g\,(g\geq 7)$ are of general type,
in particular the subvarieties of $\A_g$ of codimension one. Freitag \cite{Fr6}
showed that there exists a certain bound $g_0$ such that for $g\geq g_0,$
each irreducible subvariety of $\A_g$ of codimension one is of general type.
Weissauer \cite{Wei2} proved that every irreducible divisor of $\A_g$ is of general
type for $g\geq 10.$ Moreover he proved that every subvariety of codimension
$\leq g-13$ in $\A_g$ is of general type for $g\geq 13.$ We observe that the
smallest known codimension for which there exist subvarieties of $\A_g$ for
large $g$ which are not of general type is $g-1.\ \A_1\times \A_{g-1}$ is a
subvariety of $\A_g$ of codimension $g-1$ which is not of general type.
\vspace{0.1in}\\
\noindent{\bf Remark A.4.} Let $\M_g$ be the coarse moduli space of curves of genus $g$ over $\BC.$ Then $\M_g$ is an analytic subvariety of $\mathcal A_g$ of dimension
$3g-3.$ It is known that $\M_g$ is rational for $g=2,4,5,6.$
In 1915 Severi proved that $\M_g$ is unirational for $g\leq 10$ (see E. Arbarello and
E. Sernesi's paper \cite{A-S} for a modern rigorous proof). The unirationality of $\M_{12}$ was proved by E. Sernesi \cite{Se} in 1981. Three years later the unirationality of $\M_{11}$ and $\M_{13}$ was proved
by M.~C. Chang and Z. Ran \cite{C-R}. So the Kodaira dimension $\kappa(\M_g)$
of $\M_g$ is $-\infty$ for $g\leq 13.$ In 1982 Harris and Mumford \cite{H-M} proved that $\M_g$ is of general type  for odd $g$ with $g\geq 25$ and $\kappa(\M_{23})\geq 0.$
J. Harris \cite{H} proved that if $g\geq 40$ and $g$ is even, $\M_g$ is of general type. In 1987 D. Eisenbud and J. Harris \cite{Eis-H} proved that $\M_g$ is of general type for all $g\geq 24$ and $\M_{23}$ has the Kodaira dimension at least one. In 1996 P. Katsylo \cite{Ka} showed that $\M_3$ is rational and hence $\mathcal A_3.$
\vspace{0.1in}\\
\noindent{\bf Remark A.5.} For more details on the geometry and topology of $\mathcal A_g$ and compactifications of $\mathcal A_g$, we refer to \cite{Al1, F-C, Fr5, vdG1, vdG3, G0, G-H, I0, L-W, M-M, Sat1, Sat2, Sat3, W1}.

\vspace{0.2in}

\setcounter{theorem}{0}
\setcounter{equation}{0}
\noindent{\bf Appendix B. Extending of the Torelli Map to Toroidal Compactifications of the Siegel Modular Variety}

\vspace{0.1in}

Let ${\mathcal M}_g^{\rm DM}$ be the Deligne-Mumford compactification of $\mathcal M_g$
consisting of isomorphism classes of stable curves of genus $g$. We recall (\cite{D-M, N1, N4}) that a complete curve $C$ is said to be a {\sf stable curve} of genus $g\geq 1$ if
\medskip

\indent (S1) $C$ is reduced;\\
\indent (S2) $C$ has only ordinary double points as possible singularities;\\
\indent (S3) $\dim_\BC H^1(C,\mathcal O_C)=1$;\\
\indent (S4) each nonsingular rational component of $C$ meets the other components
at more than two points.

\vskip 3mm
P. Deligne and D. Mumford \cite{D-M} proved that the coarse moduli space
${\mathcal M}_g^{\rm DM}$ is an irreducible projective variety,and contains
$\mathcal M_g$ as a Zariski open subset.

\vskip 3mm
We have three standard explicit toroidal compactifications
$\mathcal A_g^{\rm VI},\ \mathcal A_g^{\rm VII}$ and $\mathcal A_g^{\rm cent}$
constructed from
\\
\indent (VI) the 1st Voronoi (or perfect) cone decomposition;\\
\indent (VII) the 2nd Voronoi cone decomposition;\\
\indent (cent) the central cone decomposition\\
respectively. We refer to \cite{M-V, S-B} for more details on the perfect cone decomposition and the 2nd Voronoi cone decomposition.
In 1973, Y. Namikawa \cite{N1} proposed a natural question if the Torelli map
\begin{equation*}
T_g: \mathcal M_g \lrt \mathcal A_g
\end{equation*}
extends to a regular map
\begin{equation*}
T_g^{\rm cent}: \mathcal M_g^{\rm DM} \lrt \mathcal A_g^{\rm cent}.
\end{equation*}
In fact, $\mathcal A_g^{\rm cent}$ is the normalization of the Igusa blow-up of
the Satake compactification $\mathcal A_g^{\rm Sat}$ along the boundary
$\partial\mathcal A_g^{\rm cent}$. In the 1970s, Mumford and Namikawa \cite{N2, N3}
showed that the Torelli map $T_g$ extends to a regular map
\begin{equation*}
T_g^{\rm VII}: \mathcal M_g^{\rm DM} \lrt \mathcal A_g^{\rm VII}.
\end{equation*}
In 2012, V. Alexeev and A. Brunyate \cite{Al-B} proved that the Torelli map $T_g$
can be extended to a regular map
\begin{equation*}
T_g^{\rm VI}: \mathcal M_g^{\rm DM} \lrt \mathcal A_g^{\rm VI}=\mathcal A_g^{\rm perf}
\end{equation*}
and that the extended Torelli map
\begin{equation*}
T_g^{\rm cent}: \mathcal M_g^{\rm DM} \lrt \mathcal A_g^{\rm cent}
\end{equation*}
is regular for $g\leq 6$ but not regular for $g\geq 9$. Furthermore they
also showed that the two compactifications $\mathcal A_g^{\rm VI}$ and
$\mathcal A_g^{\rm VII}$ are equal near the closure of the Jacobian locus $J_g$.
Almost at the same time the extended Torelli map $T_g^{\rm cent}$ is regular for
$g\leq 8$ by Alexeev and et al. \cite{Al-ET}.

\vskip 3mm
I would like to mention that K. Liu, X. Sun and S.-T. Yau \cite{Liu1, Liu2, Liu3, Liu4}
showed the goodness of the Hermitian metrics on the logarithmic tangent bundle
on $\mathcal M_g$ which are induced by the Ricci and the perturbed Ricci metrics
on $\mathcal M_g$. They also showed that the Ricci metric on $\mathcal M_g$
extends naturally to the divisor $D_g:={\mathcal M}_g^{\rm DM}\,\backslash\, \mathcal M_g$ and coincides with the Ricci metric on each component of $D_g$.

\vskip 3mm
Liu, Sun and Yau \cite{Liu2} showed that the existence of K{\"a}hler-Einstein metric on
$\mathcal M_g$ is related to the stability of the logarithmic cotangent bundle over
$\mathcal M_g^{\rm DM}$.

\vskip 2mm
Let $E$ be a holomorphic vector bundle over a complex manifold $X$ of dimension $n$.
Let $\Phi:=\Phi_X$ be a K{\"a}hler class (or form) of $X$. Then $\Phi$-degree of $E$ is defined by
\begin{equation*}
\deg (E):=\int_X c_1(E)\,\Phi^{n-1}
\end{equation*}
and the {\sf slope} of $E$ is defined to be
\begin{equation*}
\mu (E):= {{\deg (E)}\over {{\rm rank} (E)}}.
\end{equation*}
A bundle $E$ is said to be {\sf $\Phi$-stable} if for any proper coherent subsheaf
$\mathcal F \subset E$, we have
\begin{equation*}
\mu (\mathcal F) < \mu (E).
\end{equation*}
Let $U$ be a local chart of $\mathcal M_g$ near the boundary with pinching coordinates
$(t_1,\cdots,t_m,s_{m+1},\cdots,s_n)$ such that $(t_1,\cdots,t_m)$ represent the degeneration direction. Let
\begin{equation*}
F_i={{dt_i}\over {t_i}}\quad (1\leq i\leq m),\qquad F_j=ds_j \quad (m+1\leq j\leq n).
\end{equation*}
Then the logarithmic cotangent bundle $(T^*\mathcal M_g)^{\rm DM}$ is the unique extension of the cotangent bundle $T^*\mathcal M_g$ over $\mathcal M_g$ to
$\mathcal M_g^{\rm DM}$ such that on $U$ $F_1,F_2,\cdots,F_n$ is a local holomorphic frame of $(T^*\mathcal M_g)^{\rm DM}$.

\vspace{0.1in}
Liu, Sun and Yau \cite{Liu2} proved the following.

\vspace{0.1in}
\noindent
{\bf Theorem B.1.}
{\it The first Chern class $c_1\left( (T^*\mathcal M_g)^{\rm DM} \right)$ is positive and $(T^*\mathcal M_g)^{\rm DM}$ is stable with respect to
$c_1\left( (T^*\mathcal M_g)^{\rm DM} \right)$.}

\vspace{0.1in}
\noindent
{\bf Remark B.2.} We refer to \cite{B-Sa, Tr, Wo1, Wo2} for some topics related to
$\mathcal M_g$ and $\mathcal M_g^{\rm DM}$.

\newpage

\setcounter{theorem}{0}
\setcounter{equation}{0}
\noindent{\bf Appendix C. Singular Modular Forms}

\vspace{0.1in}

Let $\rho$ be a rational representation of $GL(g,\BC)$ on a finite
dimensional complex vector space $V_{\rho}.$ A holomorphic function
$f:\H_g\lrt V_{\rho}$ with values in $V_{\rho}$ is called a modular form
of type $\rho$ if it satisfies
$$f(M\cdot\tau)=\rho(C\tau+D)f(\tau)$$
for all $\begin{pmatrix} A & B\\ C & D \end{pmatrix}\in \G_g$ and $\tau\in \H_g.$
We denote by $[\G_g,\rho]$ the vector space of all modular forms of type
$\rho.$ A modular form $f\in [\G_g,\rho]$ of type $\rho$ has a Fourier
series $$f(\tau)=\sum_{T\geq 0}a(T)e^{2\pi i(T\tau)},\ \ \ \tau\in \H_g,$$
where $T$ runs over the set of all semipositive half-integral symmetric
matrices of degree $g.$ A modular form $f$ of type $\rho$ is said to
be {\it singular} if a Fourier coefficient $a(T)$ vanishes unless
$\text{det}\,(T)=0.$

\vspace{0.1in}
For $\tau=(\tau_{ij})\in\BH_g,$ we write $\tau=X+i\,Y$
with $X=(x_{ij}),\ Y=(y_{ij})$ real. We put
$${{\!\!\!\partial}\over {\partial Y}}=\,\left(\,{ {1+\delta_{ij}}\over 2}\,
{ {\!\!\!\!\partial}\over {\partial y_{ij}} }\,\right).$$
H. Maass \cite{M3} introduced the following differential operator
\begin{equation*}
M_g:=\det (Y)\cdot \det \left( {{\!\!\!\partial}\over {\partial Y}} \right) \leqno
({\rm C}.1)
\end{equation*}
characterizing singular modular forms. Using the differential operator $M_g$, Maass
\cite[pp.\,202--204]{M3} proved that if a nonzero singular modular form of degree $g$ and type $\rho:=\det^k$~(or weight $k$) exists, then $gk\equiv 0~({\rm mod}~2)$ and $0< 2k \leq g-1.$ The converse was proved by R. Weissauer \cite{Wei1}.
\vspace{0.1in}\\
\indent Freitag \cite{Fr3} proved that every singular modular form can be written as
a finite linear combination of theta series with harmonic coefficients and
proposed the problem to characterize singular modular forms. Weissauer \cite{Wei1}
gave the following criterion.
\vspace{0.1in}\\
\noindent{\bf Theorem C.1.}\quad{\it Let $\rho$ be an irreducible rational representation of
$GL(g,\BC)$ with its highest weight $(\l_1,\cdots,\l_g).$ Let $f$ be a
modular form of type $\rho.$ Then the following are equivalent\,:
\vspace{0.05in}\\
{\rm (a)} $f$ is singular.\vspace{0.05in}\\
{\rm (b)} $2\l_g < g.$}
\vspace{0.1in}\\
\indent Now we describe how the concept of singular modular forms is closely related
to the geometry of the Siegel modular variety. Let $X:=\mathcal A_g^{\rm Sat}$ be the Satake compactification of the Siegel modular variety $\A_g=\G_g\ba \H_g.$ Then
$\A_g$ is embedded in $X$ as a quasiprojective algebraic subvariety of
codimension $g$. Let $X_s$ be the smooth part of $\A_g$ and ${\tilde X}$ the
desingularization of $X.$ Without loss of generality, we assume
$X_s\subset {\tilde X}.$ Let $\Omega^p({\tilde X})$\,(resp.\,$\Omega^p(X_s))$
be the space of holomorphic $p$-form on ${\tilde X}$\,(resp.\,$X_s$).
Freitag and Pommerening \cite{Fr-Po} showed that if $g>1$, then the restriction
map
$$\Omega^p({\tilde X})\lrt \Omega^p(X_s)$$
is an isomorphism for $p< \text{dim}_{\BC}\,{\tilde X}={{g(g+1)}\over 2}.$
Since the singular part of $\A_g$ is at least codimension $2$ for $g>1,$
we have an isomorphism
$$\Omega^p({\tilde X})\cong \Omega^p(\H_g)^{\G_g}.$$
Here $\Omega^p(\H_g)^{\G_g}$ denotes the space of $\G_g$-invariant holomorphic
$p$-forms on $\H_g.$ Let $\text{Sym}^2(\BC^g)$ be the symmetric power of
the canonical representation of $GL(g,\BC)$ on $\BC^n.$ Then we have an
isomorphism
\begin{equation*}
\Omega^p(\H_g)^{\G_g}\lrt \left[\G_g,\bigwedge^p \text{Sym}^2(\BC^g)\right].\leqno
({\rm C}.2)
\end{equation*}

\noindent{\bf Theorem C.2.~\cite{Wei1}}
{\rm Let $\rho_{\al}$ be the irreducible representation of $GL(g,\BC)$ with highest weight
$$(g+1,\cdots,g+1,g-\al,\cdots,g-\al)$$
such that $\text{corank}(\rho_{\al})=\al$ for $1\leq \al\leq g.$ If $\al=-1,$
we let $\rho_{\al}=(g+1,\cdots,g+1).$ Then

$$\Omega^p(\H_g)^{\G_g}=\begin{cases} [\G_g,\rho_{\al}],\ &\text{if\
$p={{g(g+1)}\over 2}-{{\al(\al+1)}\over 2}$}\\
0, &\text{otherwise}.\end{cases}$$}

\noindent
{\bf Remark C.3.}
{\rm If $2\al>g,$ then any $f\in [\G_g,\rho_{\al}]$ is singular. Thus if $p<{{g(3g+2)}\over 8},$ then any $\G_g$-invariant holomorphic
$p$-form on $\H_g$ can be expressed in terms of vector valued theta series
with harmonic coefficients. It can be shown with a suitable modification
that the just mentioned statement holds for a sufficiently small congruence
subgroup of $\G_g.$}

\vspace{0.1in}
\indent Thus the natural question is to ask how to determine the $\G_g$-invariant
holomorphic $p$-forms on $\H_g$ for the nonsingular range
$\dfrac{g(3g+2)}{8}\leq p \leq \dfrac{g(g+1)}{2}.$
Weissauer \cite{Wei3} answered the above question for $g=2.$ For $g>2,$ the above question is still open. It is well know that the vector space of vector valued modular forms of type $\rho$ is finite dimensional. The computation or the estimate of the
dimension of $\Omega^p(\H_g)^{\G_g}$ is interesting because its dimension is
finite even though the quotient space $\A_g$ is noncompact.
\vspace{0.1in}\\
\indent Finally we will mention the results due to Weisauer \cite{Wei2}.
We let $\G$ be a congruence subgroup of $\G_2.$ According to
Theorem C.2, $\G$-invariant holomorphic forms in $\Omega^2(\H_2)^{\G}$ are
corresponded to modular forms of type (3,1). We note that these invariant
holomorphic $2$-forms are contained in the {\it nonsingular\ range}.
And if these modular forms are not cusp forms, they are mapped under the
Siegel $\Phi$-operator to cusp forms of weight $3$ with respect to some
congruence subgroup\,(\,dependent on $\G$\,) of the elliptic modular group.
Since there are finitely many cusps, it is easy to deal with these modular
forms in the adelic version. Observing these facts, he showed that any
$2$-holomorphic form on $\G\ba \H_2$ can be expressed in terms of theta
series with harmonic coefficients associated to binary positive definite
quadratic forms. Moreover he showed that $H^2(\G\ba \H_2,\BC)$ has a
pure Hodge structure and that the Tate conjecture holds for a suitable
compactification of $\G\ba \H_2.$  If $g\geq 3,$ for a congruence subgroup
$\G$ of $\G_g$ it is difficult to compute the cohomology groups
$H^{\ast}(\G\ba \H_g,\BC)$ because $\G\ba \H_g$ is noncompact and highly
singular. Therefore in order to study their structure, it is natural to
ask if they have pure Hodge structures or mixed Hodge structures.

\vspace{0.2in}

\setcounter{theorem}{0}
\setcounter{equation}{0}
\noindent{\bf Appendix D. Singular Jacobi Forms}

\vspace{0.1in}

In this section, we discuss the notion of singular Jacobi forms. First of all we define
the concept of Jacobi forms.

\vskip 2mm
Let $\rho$ be a rational representation of
$GL(g,\mathbb{C})$ on a finite
dimensional complex vector space
$V_{\rho}.$ Let ${\mathcal M}\in \mathbb R^{(h,h)}$ be a
symmetric half-integral semi-positive definite matrix of degree $m$.
The canonical automorphic factor
$$ J_{\rho,\mathcal M}: G^J\times \BH_{g,h}\lrt GL(V_\rho)$$
for $G^J$ on $\BH_{g,h}$ is given as follows\,:
\begin{eqnarray*}
 J_{\rho,\mathcal M}((g,(\lambda,\mu;\kappa)),(\tau,z))
&=&e^{2\,\pi\, i\,{{\rm tr}} \left( {\mathcal M}(z+\lambda\,
\tau+\mu)(C\tau+D)^{-1}C\,{}^t(z+\lambda\,\tau\,+\,\mu)\right) }\\
& &\times\, e^{-2\,\pi\, i\,{{\rm tr}} \left( {\mathcal M}(\lambda\,
\tau\,{}^t\!\lambda\,+\,2\,\lambda\,{}^t\!z+\,\kappa+
\mu\,{}^t\!\lambda) \right)} \rho(C\,\tau+D),
\end{eqnarray*}
where $g=\left(\begin{matrix} A&B\\ C&D\end{matrix}\right)\in
Sp(2g,\mathbb R),\ (\lambda,\mu;\kappa)\in H_{\mathbb R}^{(g,h)}$
and $(\tau,z)\in \BH_{g,h}.$ We refer to \cite{YJH6} for a geometrical construction of $J_{\rho,\mathcal M}.$

\vskip 0.21cm
Let
$C^{\infty}(\BH_{g,h},V_{\rho})$ be the algebra of all
$C^{\infty}$ functions on $\BH_{g,h}$
with values in $V_{\rho}.$
For $f\in C^{\infty}(\BH_{g,h}, V_{\rho}),$ we define
\begin{eqnarray*}
 \left(f|_{\rho,{\mathcal M}}[(g,(\lambda,\mu;\kappa))]\right)(\tau,z)
&= & J_{\rho,\mathcal M}((g,(\lambda,\mu;\kappa)),(\tau,z))^{-1}\\
& &\ f\left(g\!\cdot\!\tau,(z+\lambda\,
\tau+\mu)(C\,\tau+D)^{-1}\right),\nonumber
\end{eqnarray*}
where $g=\left(\begin{matrix} A&B\\ C&D\end{matrix}\right)\in
Sp(2g,\mathbb R),\ (\lambda,\mu;\kappa)\in H_{\mathbb R}^{(g,h)}$
and $(\tau,z)\in \BH_{g,h}.$

\vspace{0.1in}
\noindent
{\bf Definition D.1.}
{\rm Let $\rho$ and $\mathcal M$ be as above. Let
$$H_{\mathbb Z}^{(g,h)}:= \left\{ (\lambda,\mu;\kappa)\in H_{\mathbb R}^{(g,h)}\, \vert
\,\ \lambda,\mu, \kappa\ {\rm integral} \right\}$$
be the discrete subgroup of $H_{\mathbb R}^{(g,h)}$.
A {\sf Jacobi form} of index $\mathcal M$
with respect to $\rho$ on a subgroup $\Gamma$ of $\Gamma_g$ of finite
index is a holomorphic function $f\in C^{\infty}(\BH_{g,h},V_{\rho})$ satisfying the following conditions (A) and (B):

\smallskip

(A) \,\ $f|_{\rho,{\mathcal M}}[\tilde{\gamma}] = f$ for
all $\tilde{\gamma}\in {\widetilde\Gamma}:= \Gamma \ltimes
H_{\mathbb Z}^{(g,h)}$.

\smallskip

(B) \,\ For each $M\in\Gamma_g$, $f|_{\rho,\CM}[M]$ has a
Fourier expansion of \\
\indent \ \ \ \ \ \ \ the following form :
$$(f|_{\rho,\CM}[M])(\tau,z) = \sum\limits_{T=\,{}^tT\ge0\atop \text {half-integral}}
\sum\limits_{R\in \mathbb Z^{(g,h)}} c(T,R)\cdot e^{{ {2\pi i}\over {\lambda_\G}}
\,{{\rm tr}}(T\tau)}\cdot e^{2\pi i\,{{\rm tr}}(Rz)}$$
\indent \ \ \ \ \ \ \ with $\lambda_\G(\neq 0)\in\BZ$ and
$c(T,R)\ne 0$ only if $\left(\begin{matrix} { 1\over {\lambda_\G}}T & \frac 12 R\\
\frac 12\,^t\!R&{\mathcal M}\end{matrix}\right) \geq 0$. }

\medskip

\indent If $g\geq 2,$ the condition (B) is superfluous by K{\"o}cher 
principle\,(cf.~\cite[Lemma 1.6]{Zi}). We denote by
$J_{\rho,\mathcal M}(\Gamma)$ the vector space of all Jacobi forms
of index $\mathcal{M}$ with respect to $\rho$ on $\Gamma$.
Ziegler (cf.~\cite[Theorem 1.1]{E-Z} or \cite[Theorem 1.8]{Zi})
proves that the vector space $J_{\rho,\mathcal {M}}(\Gamma)$ is
finite dimensional. In the special case $\rho(A)=(\det(A))^k$ with
$A\in GL(g,\BC)$ and a fixed $k\in\BZ$, we write $J_{k,\CM}(\G)$
instead of $J_{\rho,\CM}(\G)$ and call $k$ the {\it weight} of the
corresponding Jacobi forms. For more results about Jacobi forms with
$g>1$ and $h>1$, we refer to \cite{YJH1, YJH4, YJH2, YJH6, YJH10, Zi}.
Jacobi forms play an important role in lifting elliptic cusp forms to 
Siegel cusp forms of degree $2g$.

\vskip 0.21cm

Without loss of generality we may assume that $\M$ is positive definite. For
simplicity, we consider the case that  $\G$ is the Siegel modular
group $\G_g$ of degree $g.$
\vspace{0.1in}\\
\indent Let $g$ and $h$ be two positive integers. We recall that $\M$ is a symmetric
positive definite, half-integral matrix of degree $h$.
We let
$${\mathcal{P}}_g:=\{ Y\in \BR^{(g,g)}\,\vert\, Y=\,{^tY} > 0 \} $$
be the open convex cone of positive definite matrices of degree $g$
in the Euclidean space $\BR^{{g(g+1)}\over 2}.$
We define the differential
operator $M_{g,h,\M}$ on ${\mathcal{P}}_g\times \BR^{(h,g)}$ defined by
\begin{equation*}
M_{g,h,{\mathcal{M}}}:=\text{det}
\,(Y)\cdot \text{det}\left( {{\partial}\over {\partial Y}}
+{1\over {8\pi}} ^t\!\left( {{\partial}\over {\partial V}}\right)
{\mathcal{M}}^{-1}
\left( {{\partial}\over
{\partial V}}\right) \right), \leqno({\rm D}.1)
\end{equation*}
where
$$Y=(y_{\mu\nu})\in {\mathcal{P}}_g,\ \ V=(v_{kl})\in \BR^{(h,g)},\ \
{{\partial}\over {\partial Y}}=\left( {{1+\delta_{\mu\nu}}\over
2}{{\partial}\over {\partial y_{\mu\nu}}}\right)$$ and
$${{\partial}\over {\partial V}}=\left( { {\partial}\over {\partial v_{kl}} }
\right).$$
We note that this differential operator $M_{g,h,{\mathcal{M}}}$ generalizes the Maass operator $M_g$~(see Formula (C.1)).
\vskip 3mm
The author \cite{YJH5} characterized singular Jacobi forms as follows\,:
\vspace{0.1in}\\
\noindent{\bf Theorem D.2.}
{\it Let $f\in J_{\rho,\M}(\G_g)$ be a Jacobi form of index
$\M$ with respect to a finite dimensional rational representation $\rho$ of
$GL(g,\BC).$ Then the following conditions are equivalent :\\
\indent {\rm (1)} $f$ is a {\it singular} Jacobi form.\\
\indent {\rm (2)} $f$ satisfies the differential equation $M_{g,h,\M}f=0.$}
\vspace{0.1in}\\
\noindent{\bf Theorem D.3.} {\it Let $\rho$ be an irreducible finite dimensional
representation of $GL(g,\BC).$ Then there exists a nonvanishing
{\it singular} Jacobi form in $J_{\rho,\M}(\G_g)$ if and only if
$2k(\rho)< g+h.$ Here $k(\rho)$ denotes the weight of $\rho.$}
\vspace{0.1in}\\
\indent For the proofs of the above theorems we refer to Theorems
4.1 and 4.5 in \cite{YJH5}.

\vspace{0.1in}
\noindent {\bf Exercise D.4.}
Compute the eigenfunctions and the eigenvalues of
$M_{g,h,\M}$ (cf.~\cite[pp.\,2048--2049]{YJH5}).
\vspace{0.1in}\\
\indent Now we consider the following group $GL(g,\BR)\ltimes
H_{\BR}^{(g,h)}$ equipped with the multiplication law $$\begin{array}{ll} &\
(A,(\l,\mu,\k))\ast (B,(\l',\mu',\k')) \\ &=(AB,\,(\l
B+\l',\mu\,^t\!B^{-1}+\mu',\k+\k'+\l B\,^t\!\mu'-\mu\,^t\!B^{-1}
\,^t\!\l')),\end{array}$$ where $A,B\in GL(g,\BR)$ and
$(\l,\mu,\k),(\l',\mu',\k')\in H_{\BR}^{(g,h)}.$ We observe that
$GL(g,\BR)$ acts on $H_{\BR}^{(g,h)}$ on the right as
automorphisms. And we have the canonical action of
$GL(g,\BR)\ltimes H_{\BR}^{(g,h)}$ on ${\mathcal{P}}_g\times
\BR^{(h,g)}$ defined by
\begin{equation*}
(A,(\l,\mu,\k))\circ (Y,V):=(AY\,^t\!A,\,(V+\l Y+\mu)\,^tA),\leqno({\rm D}.2)
\end{equation*}
where $A\in GL(g,\BR), (\l,\mu,\k)\in H_{\BR}^{(g,h)}$ and $(Y,V)\in {\mathcal{P}}_g\times
\BR^{(h,g)}.$
\vspace{0.1in}\\
\noindent{\bf Lemma D.5.}
{\it The differential operator $M_{g,h,\M}$ defined by
the formula {\rm (D.1)} is invariant under the action {\rm (D.2)} of
$GL(g,\BR)\ltimes \left\{(0,\mu,0)\,\vert\ \mu\in{\BR}^{(h,g)} \right\}.$}
\vspace{0.1in}\\
\noindent{\it Proof.}\quad It follows immediately from the direct calculation.
\vspace{0.1in}\\
\indent We have the following natural questions.
\vspace{0.1in}\\
\noindent{\bf Problem D.6.}
Develop the invariant theory for the action of
$GL(g,\BR)\ltimes H_{\BR}^{(g,h)}$ on ${\mathcal{P}}_g\times \BR^{(h,g)}.$ We refer to \cite{YJH17-1, YJH20} for related topics.
\vspace{0.1in}\\
\noindent{\bf Problem D.7.}
Discuss the application of the theory of singular
Jacobi forms to the geometry of the universal abelian variety as that of
singular modular forms to the geometry of the Siegel modular
variety (see Appendix C).

\newpage

\setcounter{theorem}{0}
\setcounter{equation}{0}
\noindent{\bf Appendix E. Stable Jacobi Forms}

\vspace{0.1in}

Throughout this appendix  we put
$$\G_g:=Sp(2g,\BZ)\quad \text{and}\quad \G_{g,h}:=\G_g \ltimes H_\BZ^{(g,h)}.$$
For a commutative ring $R$ and an integer $m$, we denote by $S_m(R)$ the set of all $m\times m$ symmetric matrices with entries in $R$.

\vskip 2mm
We know that the Siegel-Jacobi space
$$\BH_{g,h}=G^J/K^J$$
is a non-symmetric homogeneous space. Here
\begin{equation*}
K^J=\left\{ (k,(0,0;\kappa))~\vert\ k\in U(g),\ ~\kappa\in S_h(\BR) \right\}
\end{equation*}
is a subgroup of $G^J.$
Let $\frak g^J$ be the Lie algebra of the Jacobi group $G^J$.  Then
$\frak g^J$ has a decomposition
\begin{equation*}
\frak g^J =\frak k^J + \frak p^J,
\end{equation*}
where
\begin{equation*}
{\frak k}^J=\left\{\left( \begin{pmatrix} \ a & b \\ -b & a
\end{pmatrix},(0,0;\k)
\right)\ \bigg|\
a+\,{}^ta=0,\ b\in S_g(\BR),\ \k\in S_h(\BR) \right\}
\end{equation*}
and
\begin{equation*}
{\frak p}^J=\left\{\left( \begin{pmatrix} a & \ b \\ b & -a
\end{pmatrix},(P,Q;0)\right)\ \bigg|\
a,b\in S_g(\BR),\ \ P,Q\in \BC^{(h,g)} \right\}.
\end{equation*}
We observe that $\frak k^J$ is the Lie algebra of $K^J$.
The complexification ${\frak p}_\BC^J:=\frak p\otimes_\BR \BC$ of ${\frak p}^J$ has a decomposition
\begin{equation*}
\frak p^J_\BC =\frak p_+^J + \frak p_-^J,
\end{equation*}
where
\begin{equation*}
{\frak p}_+^J=\left\{\left( \begin{pmatrix} X & -iX \\ -iX & -X
\end{pmatrix},(P,-iP;0)
\right)\ \bigg|
\ X\in S_g(\BC),\ \ P\in \BC^{(h,g)} \right\}.
\end{equation*}
and
\begin{equation*}
{\frak p}_-^J=\left\{\left( \begin{pmatrix} X & -iX \\ -iX & -X
\end{pmatrix},(P,-iP;0)
\right)\ \bigg|
\ X\in S_g(\BC),\ \ P\in \BC^{(h,g)} \right\}.
\end{equation*}
We define a complex structure $I^J$ on the tangent space ${\frak p}^J$ of
$\BH_{g,h}$ at $(iI_g,0)$ by
$$I^J \left( \begin{pmatrix} a & \ b\\ b &-a
\end{pmatrix},(P,Q;0)\right):=
\left(\begin{pmatrix} \ b & -a\\ -a & -b
\end{pmatrix},(Q,-P;0)\right).$$
Identifying $\BR^{(h,g)}\times \BR^{(h,g)}$ with $\BC^{(h,g)}$ via
$$(P,Q;0)\lmt i\,P+Q,\ \ P,Q\in \BR^{(h,g)},$$
we may regard the complex structure $I^J$ as a real linear map
$$I^J(X+i~Y,Q+i~P)=(-Y+i\,X,-P+i\,Q),$$
where $X+i~Y\in S_g(\BR),\ Q+i\,P\in \BC^{(h,g)}.\ I^J$
extends complex linearly on the complexification ${\frak p}_{\BC}^J.$
With respect to this complex structure $I^J,$ we may say that a function $f$
on $\BH_{g,h}$ is {\it holomorphic} if and only if $\xi f=0$ for all
$\xi\in {\frak p}_-^J.$

\vskip 3mm
Since the space $\BH_{g,h}$ is diffeomorphic to the homogeneous space
$G^J/K^J$, we may lift a function $f$ on $\BH_{g,h}$ with values in $V_{\rho}$
to a function $\Phi_f$ on $G^J$ with values in $V_{\rho}$ in the following
way. We define the lifting
$$L_{\rho,\M}\,:\,{\mathcal F}(\BH_{g,h},V_{\rho})\lrt
{\mathcal F}(G^J,V_{\rho}),\ \ \ L_{\rho,\M}(f):=\Phi_f \leqno({\text E}.1)$$
by
\begin{eqnarray*}
\Phi_f(x):\!\!\!\!&=&\!\!\!\! (f|_{\rho,\M}[x])(iI_g,0)\\
            &=&\!\!\!\! J_{\rho,\M}(x,(iI_g,0))\,f(x\cdot (iI_g,0)),
\end{eqnarray*}
where $x\in G^J$ and ${\mathcal F}(\BH_{g,h},V_{\rho})\,(\text{resp.}
\ {\mathcal F}(G^J,V_{\rho}))$ denotes the vector space consisting of functions on $\BH_{g,h}\,(\text{resp.}\ G^J)$ with values in $V_{\rho}$.

\vskip 3mm
We see easily that the vector space $J_{\rho,\M}(\G_g)$ is
isomorphic to the space $A_{\rho,\M}(\G_{g,h})$ of smooth functions $\Phi$ on
$G^J$ with values in $V_{\rho}$ satisfying the following conditions:
\vskip 2mm
(1a)\ \ \ $\Phi(\g x)=\Phi(x)$\ \ for\ all\ $\g\in \G^J$ and $x\in G^J.$\par
(1b)\ \ \ $\Phi(x\,r(k,\kappa))=e^{2\pi i\,\s (\M\kappa)}\rho(k)^{-1}
\Phi(x)$ \ \ \ for\ all\ $x\in G^J,\ r(k,\kappa)\in K^J.$\par
(2)\ \ \ $Y^{-}\Phi=0$\ \ \ for all $Y^{-}\in {\frak p}^J_-$.\par
(3)\ \ For all $M\in Sp(2g,\BR),$ the function $\psi:G^J\lrt V_{\rho}$
defined by
$$\psi(x):=\rho(Y^{-{1\over 2}})\,\Phi(Mx),\ \ \ x\in G^J$$
\indent \ \ \ \ \ \,
is bounded in the domain $Y\geq Y_0.$ Here $x\cdot (iI_g,0)=(\tau,z)$ with
$\tau=$\par
\indent \ \ \ \ \ \, $X+i\,Y,\ Y>0.$ 

\vskip 3mm
Clearly $J_{\rho,\M}^{\text{cusp}}(\G_g)$
is isomorphic to the
subspace $A^0_{\rho,\M}(\G_{g,h})$ of $A_{\rho,\M}(\G_{g,h})$
with the condition
(3+) that the function $g\longmapsto \Phi(g)$ is bounded.

\vskip 3mm
Let $\M$ be a fixed positive definite symmetric half-integral matrix of
degree $h$. Let $\rho_{\infty}:=\,(\rho_n)$ be a stable representation
of $GL(\infty,\BC).$ That is, for each $n\in \BZ^+,\ \rho_n$ is a finite
dimensional rational representation of $GL(n,\BC)$ and $\rho_{\infty}$
is compatible with the embeddings $\al_{kl}:GL(k,\BC)\lrt GL(l,\BC)\,
(\,k<l\,)$ defined by
$$\al_{kl}(A):=\,\begin{pmatrix} A & 0 \\ 0 & I_{l-k}\end{pmatrix},
\ \ \ A\in GL(k,\BC),\ \ k<l.$$
For two positive integers $m$ and $n$, we put
$$G_{n,m}^J:=Sp(2n,\BR)\ltimes H_\BR^{(n,m)}.$$
\def\tA{\tilde A}
For $k,l\in \BZ^+$ with $k<l,$ we define the mapping $\Phi_{l,k,\M}$
of $A_{\rho_l,\M}(\G_{l,m})$ into the functions on
$G_{k,m}^J$ by
$$\left(\,\Phi_{l,k,\M}F\,\right)(x):=\,J_{\M,\rho_k}(x,\,(iI_k,0))\,
\lim_{t\lrt\infty}\,J_{\M,\rho_l}(x_t,\,(iI_l,0))^{-1}F(x_t),
\leqno({\rm E}.2)$$
where $F\in A_{\rho_l,\M}(\G_{l,m}^J),\ x=(M,(\l,\mu;\k))\in G_{k,m}^J$
with $M=\begin{pmatrix} A & B\\ C & D\end{pmatrix}\in Sp(2k,\BR)$ and
$$x_t:=\,\left(\,\begin{pmatrix} A & 0 & B & 0\\
0 & t^{1/2}I_{l-k} & 0 & 0 \\ C & 0 & D & 0 \\
0 & 0 & 0 & t^{-1/2}I_{l-k}\end{pmatrix},\,(\,(\l,0),(\mu,0);\k)\right)
\in G_{l,h}^J.$$
\noindent
{\bf Proposition\ {\rm E}.1.} {\it The limit {\rm (E.2)} always exists and the image
of $A_{\rho_l,\M}(\G_{l,h})$ under $\Phi_{l,k,\M}$ is contained in
$A_{\rho_k,\M}(\G_{k,h})$. Obviously the mapping
$$\Phi_{l,k,\M}:A_{\rho_l,\M}(\G_{l,h}^J)\lrt A_{\rho_k,\M}(\G_{k,h})$$
is a linear mapping.}
\vskip 2mm
The mapping $\Phi_{l,k,\M}$ is called the {\sf Siegel-Jacobi operator}.
For any $g\in \BZ^+,$ we put
$$A_{g,\M}:=\,\bigoplus_{\rho} A_{\rho,\M}(\G_{g,h}),\leqno ({\rm E}.3)$$
where $\rho$ runs over all isomorphism classes of irreducible rational
representations of $GL(g,\BC)$. For $g=0,$ we set $A_{0,\M}:=\BC.$
\vskip 2mm
For each $g\in \BZ^+,$ we put
$$A_{g,\M}^{\ast}:=\,\bigoplus_{\rho_{\ast}} A(\rho_{\ast},\M),\leqno ({\rm E}.4)$$
where $\rho_{\ast}$ runs over all isomorphism classes of irreducible
rational representations of $GL(g,\BC)$ with highest weight
$\lambda(\rho_{\ast})\in (2\BZ)^g.$ It is obvious that if $k<l,$ then the
Siegel-Jacobi operator $\Phi_{l,k,\M}$ maps $A_{l,\M}\,$(\,resp.\,
$A_{l,\M}^{\ast}\,)$ into $A_{k,\M}$\,(\,resp.\,$A_{k,\M}^{\ast}\,).$

\vskip 2mm
We let
$$A_{\infty,\M}:=\,\lim_{\begin{subarray}{c} \longleftarrow\\ ^k \end{subarray}}
A_{k,\M}\ \ \
\text{and}\ \ \ A_{\infty,\M}^{\ast}:=\,
\lim_{\begin{subarray}{c} \longleftarrow\\ ^k \end{subarray}}
A_{k,\M}^{\ast}  \leqno ({\rm E}.5)$$
be the inverse limits of $(\,A_{k,\M},\,\Phi_{l,k,\M}\,)$ and
$(\,A_{k,\M}^{\ast},\,\Phi_{l,k,\M}\,)$ respectively.

\vskip 2mm
\noindent
{\bf Proposition\ {\rm E}.2.} {\it $A_{\infty,\M}$ has a commutative ring
structure compatible with the Siegel-Jacobi operators. Obviously
$A_{\infty,\M}^{\ast}$ is a subring of $A_{\infty,\M}.$}

\vskip 2mm
For a stable irreducible representation $\rho_{\infty}=(\rho_g)$ of
$GL(\infty,\BC)$, we define
$$A_{\rho_{\infty},\M}:=
\,\lim_{\begin{subarray}{c} \longleftarrow\\ ^g \end{subarray}}
A_{\rho_g,\M}(\G_{g,h}).  \leqno ({\rm E}.6)$$
\noindent
{\bf Proposition\ {\rm E}.3.} {\it We have
$$A_{\infty,\M}\,=\,\bigoplus_{\rho_{\infty}}\,A_{\rho_{\infty},\M},$$
where $\rho_{\infty}$ runs over all isomorphism classes of stable
irreducible representations of $GL(\infty,\BC).$}

\vskip 3mm
\noindent
{\bf Definition\ {\rm E}.4.} {\rm Elements in $A_{\infty,\M}$ are called
{\sf stable automorphic forms} on $G_{\infty,h}^J$ of index
$\M$ and elements of $A_{\infty,\M}^{\ast}$ are called
{\sf even stable automorphic forms} on $G_{\infty,h}^J$ of
index $\M$.}

\vskip 2mm
For $g\geq 1,$ we define
$${\mathbb A}_g:=\,\bigoplus_{\rho}\bigoplus_{\M} A_{\rho,\M}(\G_{g,h}),
\leqno ({\rm E}.7)$$
where $\rho$ runs over all isomorphism classes of irreducible
rational representations of $GL(g,\BC)$ and $\M$ runs over all
equivalence classes of positive definite symmetric, half-integral
matrices of any degree $\geq 1.$ We set ${\mathbb A_0}:=\BC.$

\vskip 2mm
For $g\geq 1,$ we also define
$${\mathbb A}_g^{\ast}:=\,\bigoplus_{\rho_{\ast}}
\bigoplus_{\M} {\mathbb A}_{\rho_{\ast},\M}(\G_{g,h}),
\leqno ({\rm E}.8)$$
where $\rho_{\ast}$ runs over all isomorphism classes of irreducible
rational representations of $GL(g,\BC)$ with highest weight
$\lambda(\rho_{\ast})\in (2\BZ)^g$ and $\M$ runs over all equivalence classes of positive definite symmetric half-integral matrices of any degree $\geq 1.$

\vskip 2mm
Let $\rho_{\infty}=(\rho_g)$ be a stable irreducible rational
representation of $GL(\infty,\BC).$ For each irreducible rational
representation $\rho_g$ of $GL(g,\BC)$ appearing in $\rho_{\infty},$
we put
$$A(\rho_g;\rho_{\infty}):=\,\bigoplus_{\M}
A_{\rho_g,\M}(\G_{g,h}),\leqno ({\rm E}.9)$$
where $\M$ runs over all equivalence classes of positive definite
symmetric half-integral matrices of any degree $\geq 1.$ Clearly
the Siegel-Jacobi operator $\Phi_{l,k}:=\,\bigoplus_{\M}
\Phi_{l,k,\M}\,(\,k< l\,)$ maps $A(\rho_l;\rho_{\infty})$ into
$A(\rho_k;\rho_{\infty}).$

\vskip 2mm
Using the Siegel-Jacobi operators, we can define the inverse limits
$$A(\rho_{\infty}):=\,\lim_{\begin{subarray}{c} \longleftarrow\\ ^g \end{subarray}}
A(\rho_g;\rho_{\infty}),
\ \ \ A_{\infty}:=\,\lim_{\begin{subarray}{c} \longleftarrow\\ ^g \end{subarray}}
A_g\ \ \ \text{and}\ \
A_{\infty}^{\ast}:=\,\lim_{\begin{subarray}{c} \longleftarrow\\ ^g \end{subarray}}
A_g^{\ast}.\leqno ({\rm E}.10)$$
\noindent
{\bf Theorem\ {\rm E}.5.}
{\it $$A_{\infty}\,=\,\bigoplus_{\rho_{\infty}} A(\rho_{\infty}),$$
where $\rho_{\infty}$ runs over all equivalence classes of stable
irreducible representations of $GL(\infty,\BC).$}

\vskip 3mm

\indent  Let $\rho$ and $\M$ be the same as in the previous sections. For positive
integers $r$ and $g$ with $r<g,$ we let $\rho^{(r)}:GL(r,\BC)\lrt
GL(V_{\rho})$ be a rational representation of $GL(r,\BC)$ defined by
$$\rho^{(r)}(a)v:=\rho\left( \begin{pmatrix} a & 0\\ 0 & I_{g-r}\end{pmatrix}
\right)v,\ \ \ a\in GL(r,\BC),\ \ v\in V_{\rho}.$$
The Siegel-Jacobi operator $\Psi_{g,r}:J_{\rho,\M}(\G_g)\lrt
J_{\rho^{(r)},\M}(\G_r)$ is defined by
\begin{equation*}
(\Psi_{g,r}f)(\tau,z):=\lim_{t\rightarrow \infty}\,f\left(
\begin{pmatrix} \tau & 0\\ 0 & itI_{g-r}\end{pmatrix},(z,0)\right),\leqno ({\rm E}.11)
\end{equation*}
where $f\in J_{\rho,\M}(\G_g),\ \tau\in \H_r$ and $z\in \BC^{(h,r)}.$
It is easy to check that the above limit always exists and the
Siegel-Jacobi operator is a linear mapping. Let $V_{\rho}^{(r)}$ be the
subspace of $V_{\rho}$ spanned by the values
$\{\,(\Psi_{g,r}f)(\tau,z)\,\vert\ f\in J_{\rho,\M}(\G_g),\ (\tau,z)\in
\H_r\times \BC^{(h,r)}\,\}.$ Then $V_{\rho}^{(r)}$ is invariant under the
action of the group
$$\left\{\,\begin{pmatrix} a & 0\\ 0 & I_{g-r}\end{pmatrix}\,:\ a\in GL(r,\BC)\
\right\}\cong GL(r,\BC).$$
We can show that if $V_{\rho}^{(r)}\neq 0$ and $(\rho,V_{\rho})$ is
irreducible, then $(\rho^{(r)},V_{\rho}^{(r)})$ is also irreducible.
\vspace{0.1in}\\
\noindent{\bf Theorem E.6.}
{\it The action of the Siegel-Jacobi operator is compatible
with that of that of the Hecke operator.}
\vspace{0.1in}\\
\indent We refer to \cite{YJH4} for a precise detail on the Hecke operators and
the proof of the above theorem.
\vspace{0.1in}\\
\noindent{\bf Problem E.7.} {\rm Discuss the injectivity, surjectivity and bijectivity
of the Siegel-Jacobi operator.}
\vspace{0.1in}\\
\indent This problem was partially discussed by the author \cite{YJH4} and Kramer \cite{K} in the
special cases. For instance, Kramer \cite{K} showed that if $g$ is arbitrary, $h=1$ and $\rho:GL(g,\BC)\lrt \BC^{\times}$ is a one-dimensional representation of $GL(g,\BC)$
defined by $\rho(a):=(\text{det}\,(a))^k$ for some $k\in \BZ^+,$ then
the Siegel-Jacobi operator $$\Psi_{g,g-1}:J_{k,m}(\G_g)\lrt J_{k,m}(\G_{g-1})$$
is surjective for $k\gg m\gg 0.$
\vspace{0.1in}\\
\noindent{\bf Theorem E.8.}
{\it Let $1\leq r\leq g-1$ and let $\rho$ be an irreducible
finite dimensional representation of $GL(g,\BC).$ Assume that
$k(\rho)> g+r+\text{rank}\,(\M)+1$ and that $k$ is even. Then
$$J_{\rho^{(r)},\M}^{\text{cusp}}(\G_r)\subset
\Psi_{g,r}(J_{\rho,\M}(\G_g)).$$
Here $J_{\rho^{(r)},\M}^{\text{cusp}}(\G_r)$ denotes the subspace consisting
of all cuspidal Jacobi forms in $J_{\rho^{(r)},\M}(\G_r).$}
\vspace{0.1in}\\
\noindent{\it Idea of Proof.}\quad For each $f\in J_{\rho^{(r)},\M}^{\text{cusp}}
(\G_r),$ we can show by a direct computation that
$$\Psi_{g,r}(E_{\rho,\M}^{(g)}(\tau,z;f))=f,$$
where $E_{\rho,\M}^{(g)}(\tau,z;f)$ is the Eisenstein series of Klingen's
type associated with a cusp form $f.$ For a precise detail, we refer to \cite{Zi}.
\vspace{0.1in}\\
\noindent{\bf Remark E.9.}
Dulinski \cite{Du} decomposed the vector space
$J_{k,\M}(\G_g)\,(k\in \BZ^+)$ into a direct sum of certain subspaces
by calculating the action of the Siegel-Jacobi operator on
Eisenstein series of Klingen's type explicitly.
\vspace{0.1in}\\
\indent For two positive integers $r$ and $g$ with $r\leq g-1,$ we consider the
bigraded ring
$$J_{\ast,\ast}^{(r)}(\ell):=\bigoplus_{k=0}^{\infty}\bigoplus_{\M}
J_{k,\M}(\G_r(\ell))$$
and
$$M_{\ast}^{(r)}(\ell):=\bigoplus_{k=0}^{\infty}\,J_{k,0}(\G_r(\ell))
=\bigoplus_{k=0}^{\infty}[\G_r(\ell),k],$$
where $\G_r(\ell)$ denotes the principal congruence subgroup of $\G_r$
of level $\ell$ and $\M$ runs over the set of all symmetric semi-positive
half-integral matrices of degree $h$.
Let
$$\Psi_{r,r-1,\ell}:J_{k,\M}(\G_r(\ell))\lrt J_{k,\M}(\G_{r-1}(\ell))$$
be the Siegel-Jacobi operator defined by (E.11).
\vspace{0.1in}\\
\noindent{\bf Problem E.10.}\quad Investigate $\text{Proj}\,J_{\ast,\ast}^{(r)}(\ell)$ over
$M_{\ast}^{(r)}(\ell)$ and the quotient space
$$Y_r(\ell):=
({{\G_r(\ell)\ltimes (\ell\BZ)^2})\ba ({H_r\times \BC^r}})$$ for
$1\leq r\leq g-1.$\vspace{0.1in}\\
The difficulty to this problem comes from the following facts (A) and
(B)\,:
\vspace{0.1in}

\indent (A) $J_{\ast,\ast}^{(r)}(\ell)$ is not finitely generated over
$M_{\ast}^{(r)}(\ell).$\\
\indent (B) $J_{k,\M}^{\text{cusp}}(\G_r(\ell))\neq \text{ker}\,\Psi_{r,r-1,\ell}$
in general.
\vspace{0.1in}\\
\indent These are the facts different from the theory of Siegel modular forms.
We remark that Runge~(cf.~\cite[pp.\,190--194]{Ru}) discussed some parts about the above problem.

%
%
%
%
%
%
%
%
%
\footnotesize{

}

\end{document}